 \documentclass[draft,3p,times]{elsarticle}

\usepackage{amssymb}
 \usepackage{amsthm}
 \usepackage[T1]{fontenc}
\usepackage[ansinew]{inputenc}
\usepackage{mathrsfs}
\usepackage{euscript}
\usepackage{natbib}

\newtheorem{theorem}{Theorem}
\newtheorem{corollary}{Corollary}
\newtheorem{definition}{Definition}
\newtheorem{lemma}{Lemma}
\newdefinition{rmk}{ Remark}
\newproof{pf}{Proof}

\begin{document}

\begin{frontmatter}



\title{ Finite time blow-up for a wave equation with a nonlocal nonlinearity}

 \author[label1,label2,label3]{A. Z. Fino}
  \ead{ahmad.fino01@gmail.com}
 \address[label1]{Laboratoire de math\'ematiques appliqu\'ees, UMR CNRS 5142, Universit\'e de Pau et des Pays
de l'Adour, 64000 Pau, France}
\address[label2]{LaMA-Liban, Lebanese University, P.O. Box 37 Tripoli, Lebanon}
\author[label3]{M. Kirane}
\ead{mokhtar.kirane@univ-lr.fr}
\address[label3]{D\'{e}partement de Math\'{e}matiques, Universit\'{e} de la Rochelle, 17042 La Rochelle, France}
\author[label4]{V. Georgiev}
\ead{georgiev@dm.unipi.it}
\address[label4]{UNIVERSIT\'A DI PISA, Dipartimento di Matematica, "L. Tonelli" $56127$ Pisa, Largo Bruno Pontecorvo $5,$ Italy}

 \begin{abstract}
In this article, we study the local existence of solutions for a wave equation with a nonlocal in time nonlinearity. Moreover, a blow-up results are proved under some conditions on the dimensional space, the initial data and the nonlinear forcing term.
\end{abstract}


\end{frontmatter}


\section{Introduction}
\setcounter{equation}{0}

 We study the following nonlinear
wave type equation which contains a nonlocal in time nonlinearity
\begin{equation}\label{1++}
\begin{array}{ll}
\,\,\displaystyle{ u_{tt}-\Delta
u=\frac{1}{\Gamma(1- \gamma)}
\int_0^t(t-s)^{-\gamma}|u(s)|^p\,ds}&\displaystyle {x\in {\mathbb{R}^N},t>0,}\\
\end{array}
 \end{equation}
where $0<\gamma<1,$ $p>1,$ $N\geq1,$
$\Delta$ is the standard Laplacian and  $\Gamma$ is
the Euler gamma function. The nonlinear nonlocal term can be considered as an approximation of the classical semilinear
wave equation
$$u_{tt}-\Delta
u=|u(t)|^p\, $$
since the limit
$$ \lim_{\gamma \rightarrow 1} \frac{1}{\Gamma(1- \gamma)} s_+^{-\gamma} = \delta(s)$$
exists in distribution sense.

\noindent It is clear that this nonlinear term involves memory type selfinteraction and can be considered  as Riemann-Liouville integral operator
$$
_a D_t^{-\alpha} = J^\alpha_{a|t}g(t):=\frac{1}{\Gamma(\alpha)}\int_a^t(t-s)^{\alpha-1}g(s)\,ds
$$
introduced with $a=-\infty$
by Liouville in 1832 and with $a=0$
 by Riemann in 1876 (see Chapter V in \cite{DB}).
Therefore, (\ref{1++}) takes the form

\begin{equation}\label{A.2++}
     u_{tt}-\Delta
u=D_t^{\gamma - 1}\left(|u(t)|^p\right),
\end{equation}
where $ D_t^{-\alpha} =  J^\alpha_{0|t}$ and $\alpha = 1 - \gamma.$

\noindent In this work we study blow up phenomena for this semilinear wave equation and small initial  data
\begin{equation}\label{A.2++rr}
\begin{array}{ll}
\,\,
\displaystyle{u(0,x)=u_0(x),\;u_t(0,x)=u_1(x)}&\displaystyle{x\in {\mathbb{R}^N},}
\end{array}
\end{equation} where $$(u_0,u_1)\in \mathcal{H}^\mu = H^\mu(\mathbb{R}^N)\times H^{\mu-1}(\mathbb{R}^N)$$ 
and $H^\mu(\mathbb{R}^N)$ is the classical Sobolev space of order $\mu>0.$\\

The study of the non-existence of global solutions to semilinear wave equations has been initiated in the early sixties by Keller and intensively developed since then by John and Kato. It is based on an averaging method for positive solutions, usually with compact support. Much has been devoted to the case of the equation
\begin{equation}\label{simple1++}
u_{tt}-\Delta u=|u|^p,\quad p>1.
\end{equation}
It is well known that this problem does not admit a global solution for any $p>1$ when the initial values $u_0$ and $u_1$ are large in some sense (cf. \cite{Glassey,Keller,Levine}). On the other hand, John proved in \cite{John}, when $N=3,$ that nontrivial solutions with compactly supported initial data must
blow up in finite time when $1 < p < 1+\sqrt{2}.$ Interestingly, Strauss discovered the same number as the root of a dimension dependent polynomial in his work on low energy scattering for the nonlinear Klein-Gordon equation \cite{Strauss}. This led him to conjecture that the critical value, $p_0(N),$ generalizing John's result to $N$ dimensions, should be the positive root of
$$(N-1)p^2-(N+1)p-2=0.$$
Glassey \cite{Glassey} verified the conjecture when $N = 2$ under the additional assumption that $u_0$ and $u_1$ have both positive average. The technique used by Glassey, John and Sideris is to derive differential inequalities which are satisfied by the average function $t\longmapsto\int_{\mathbb{R}^N}u(x,t)\,dx.$ The fact that the support of $u(\cdotp,t)$ is included in the cone $\{x;\;|x|< t+R\}$ plays a fundamental role in deriving the differential inequalities.\\
Sideris \cite{Sideris} completes this conjecture for $N>3$ and proved that global solutions do not exist when $1<p<p_0(N),$ provided that the initial data are compactly supported and satisfy the positivity condition
$$\int_{\mathbb{R}^N}|x|^{\eta-1}u_0>0\quad\mbox{ and} \quad\int_{\mathbb{R}^N}|x|^\eta u_1>0,$$
where $\eta=0$ if $N$ is odd and $1/2$ if $N$ is even.\\
The critical case $p=p_0(N)$ was studied by Schaffer \cite{Schaeffer} in dimension $N=2$ and $N=3,$ and then completed in $2006$ by Yordanov and Zhang \cite{Yordanov} for the case $N\geq4.$ \\
A slightly less sharp result under much weaker assumptions was obtained by Kato \cite{Kato} with a much easier proof. In particular, Kato pointed out the role of the exponent ${(N+1)}/{(N-1)}<p_0(N),$ for $N\geq2,$ in order to have more general initial data, but still with compact support. \\

In this paper, we generalize Kato and Glassey-Strauss critical exponents and give sufficient conditions for finite time blow-up of a new type of class of equations $(\ref{1++})$ with nonlocal in time nonlinearities. Let us mention that our blow-up results and initial conditions are similar to that of Kato and Glassey-Strauss respectively.\\

Our first point to discuss the existence of local solutions to $(\ref{1++})$ with initial data $(\ref{A.2++rr})$.
Formally, the equation $(\ref{simple1++})$ can be rewritten as integral equation
$$ u(t)=\dot{K}(t)u_0+K(t)u_1+ N(u)(t), ,\quad t\in[0,T],$$
where $K(t) = \omega^{-1}\sin\omega t, $  $\omega:=(-\Delta)^{1/2}$ and
$$ N(u) = \int_0^t K(t-s)  D_t^{\gamma - 1}(|u|^{p})(s) ds.$$
The general setting for the well -- posedness of this integral equation with
$(u_0,u_1) \in \mathcal{H}^\mu$ requires to define for any $T  >0$  a closed subspace
$$ X(T) \subseteq C([0,T],H^\mu(\mathbb{R}^N))\cap C^1([0,T],H^{\mu -1}(\mathbb{R}^N)) $$
such that
$$ (u_0,u_1) \in \mathcal{H}^\mu  \Longrightarrow  \dot{K}(t)u_0 \in X(T), \ K(t)u_1 \in X(T)$$
and
$$ u \in X(T)  \Longrightarrow  N(u) \in X(T).$$
Then the integral equation is well -- posed in $\mathcal{H}^\mu,$ if for any $R>0$ one can find $T=T(R) >0$ so that
for any initial data satisfying
$$ \|(u_0,u_1)\|_{\mathcal{H}^\mu} \leq R,$$
the integral equation
$$ u(t)=\dot{K}(t)u_0+K(t)u_1+ N(u)(t), \quad t\in[0,T],$$
has a unique solution $u \in X(T).$ Once the well posedness of the integral equation is established, one can easily prove
there exist a maximal
time $T_{\max}>0$ and a unique  solution $u \in X(T)$ for any $T \in [0,T_{max})$, such that if $T_{\max}< \infty$, then
$$ \lim_{t \nearrow T_{max}} \|u(t)\|_{H^\mu} + \| u_t(t) \|_{H^{\mu -1}} = \infty $$
i.e.  the  $\mathcal{H}^\mu$ -- norm of the solution blows up at $t=T_{max}.$

When $\mu = 1$ and $ p$ satisfies
$$\left\{
\begin{array}{ll}
\displaystyle{1 < p \leq \frac{N}{N-2}}&\displaystyle{\quad\mbox{if}\;\;N>2}\\
\displaystyle{1 < p < \infty }&\displaystyle{\quad\mbox{if}\;\;N=1,2}. \\
 \end{array}
 \right.$$
one can take $$ X(T) = C([0,T],H^1(\mathbb{R}^N))\cap C^1([0,T],L^2(\mathbb{R}^N)) $$ and using contraction mapping principle
to obtain unique solution $u \in X(T).$ (see our Theorem $\ref{Theoremmild0++}$ in Section $\ref{sec3++}$ below)
These type of solutions are called mild  solutions and the proof of the existence and uniqueness of mild solutions needs only energy type estimates and Sobolev embeddings.

The interval $ p \in (1, N/(N-2) )$ is not optimal for the local existence of solutions, but enables us to obtain first blow-up results. To state them we first define
\begin{equation}\label{conditionoverp1}
p_1=p_1(N,\gamma) := 1+\frac{3-\gamma}{(N-2+\gamma)},
\end{equation}
so that $p_1$ is the Kato exponent for $\gamma =1.$ The other quantity that generalizes Glassey-Strauss exponent (at least for $N=3$) and it  is the positive root $p_2=p_2(N,\gamma)$ of the equation
\begin{equation}\label{p_2}
(N-2)p^2-(N-\gamma)p-1=0,\quad N \geq 3.
\end{equation}

Taking $N=3$ one can see that standard observation that Kato's exponent is below the exponent of Glassey - Strauss, might be not true if $\gamma$ varies in the interval $(0,1).$  Indeed
$$ \lim_{\gamma \nearrow 1} p_1(3, \gamma) = 2 < \lim_{\gamma \nearrow 1} p_2(3, \gamma) = 1 + \sqrt{2},$$
$$  p_1(3, 1/3) = 3 = p_2(3,1/3)$$
while
$$ \lim_{\gamma \searrow 0} p_1(3, \gamma) = 4 > \lim_{\gamma \searrow 0} p_2(3, \gamma) = \frac{3 + \sqrt{13}}{2}.$$

Our first blow up result treats the case $ \gamma \in [1/3,1), $ since in this case we have
$$p_1(3,\gamma) \leq  p_2(3,\gamma) \leq N/(N-2)=3,$$
i.e. local existence requirements for mild solutions in energy space are satisfied.

Then we have the following blow up result.

\begin{theorem}\label{Blow-up-theorem2N3} Suppose $ \gamma \in [1/3,1) $ and $(u_0,u_1)\in H^1(\mathbb{R}^3)\times L^2(\mathbb{R}^3)$ satisfy $$\mbox{supp}u_i\subset B(r):=\{x\in\mathbb{R}^n:\;|x|< r\},\quad r>0,\;i=0,1, $$ and
\begin{equation}\label{averages1N3}
 \int_{\mathbb{R}^3}u_1>0,\qquad\mbox{and}\qquad \int_{\mathbb{R}^3}|x|^{-1}u_0>0,
\end{equation}
If $p< p_2=p_2(3,\gamma),$ where $p_2$ is given in $(\ref{p_2}),$ then the solution of $(\ref{1++})$
blows up in finite time.\\
\end{theorem}
Our more general result for the case $N \geq 5$ odd and $ \gamma \in [(N-2)/N,1) $
can be found in Theorem $\ref{Blow-up-theorem2}.$

Turning to the case $N=4$ we can use the following property
$$ \gamma \in [1/2,1) \Longrightarrow p_2(4,\gamma) < p_1(4,\gamma) \leq N/(N-2)= 2. $$
The corresponding blow up result reads as.

\begin{theorem}\label{blowuptheorem++N4} Assume $\gamma \in [1/2,1), N=4$  and let $(u_0,u_1)\in H^1(\mathbb{R}^4)\times L^2(\mathbb{R}^4)$ be such that
\begin{equation}\label{averages4}
\int_{\mathbb{R}^4}u_0>0,\quad \int_{\mathbb{R}^4}u_1>0.
\end{equation}
If $p\leq p_1 = p_1(4,\gamma),$ where $p_1$ is given in $(\ref{conditionoverp1}),$ then the solution of $(\ref{1++})$
blows up in finite time.\\
\end{theorem}
 The generalization of this result for the case $N \in\{1\}\cup\{2m,\;m\in\mathbb{N}^*\}$  and $ \gamma \in [(N-2)/N,1)$ is presented in Theorem $\ref{blowuptheorem++}.$

To treat values of $\gamma>0$ such that
$\gamma < {(N-2)}/{N},$ for $N\geq3,$ one has to take into account the fact that we have ${N}/{(N-2)}< p_2<p_1<1/\gamma,$
so mild solutions with data in the energy space and nonlinear exponent $ p \leq N/(N-2)$ are not sufficient to obtain blow up result for all values
of $p \in (0,1/\gamma]$ and all $\gamma \in (0, {(N-2)}/{N}).$ One slight improvement of the requirements on $p$ for the local well posedness can be done if we consider mild solutions with initial data of higher regularity, i.e.
$(u_0,u_1) \in \mathcal{H}^\mu$ with $\mu >1.$ Then the mild solution have to belong to the space
$$ X(T) = C([0,T],H^\mu(\mathbb{R}^N))\cap C^1([0,T],H^{\mu -1}(\mathbb{R}^N)). $$
When $N\geq 3$ and $ p>N/(N-2)$ satisfies
\begin{equation}\label{eq.mildden}
   \begin{array}{ll}
\displaystyle{p^2 - \frac{pN}{2}+ \frac{N}{2} \geq  0 } \\
 \end{array}
\end{equation}
one can  use contraction mapping principle
to obtain unique solution $u \in X(T)$ with $\mu = N/2-1/(p-1)$(see our Theorem $\ref{TheoremmildP0++}$ in Section $\ref{sec3++}$ below). The result is established by using only energy type estimates and Sobolev embedding. The condition
$(\ref{eq.mildden})$ is always true for space dimensions $3 \leq N \leq 8,$ but is still very restrictive for higher dimensions.

To cover larger interval for $p$ where local existence and uniqueness can be established we take
$$ X(T) = C([0,T],H^1(\mathbb{R}^N))\cap C^1([0,T],L^2(\mathbb{R}^N)) \cap
L^{q}([0,T];L^{r}_x), $$
where $(q,r)$ and the spaces $ L^{q}([0,T];L^{r}_x)$ are involved in the Stichartz estimates for the wave equation.

Note that similar spaces have been used by Ginibre and Velo in \cite{GV1} and \cite{GV2}, where
the local well posedness of the Cauchy problem for the semilinear wave equation is studied under the assumption
$ p \leq (N+2)/(N-2).$
In our case of nonlinear memory type term we are able to establish the following.

\begin{theorem}\label{Theoremweak0++in}
Given $(u_0,u_1)\in H^1(\mathbb{R}^N) \times L^2(\mathbb{R}^N),$ $N\geq 1,$  $\gamma\in(0,1)$ and let $p>1$ be such that
$$
\left\{
\begin{array}{lll}
\displaystyle{1 < p < \infty }&\displaystyle{\quad\mbox{if}\;\;N=1,2}, \\
\displaystyle{1 < p < \frac{N+4-2\gamma}{N-2}}&\displaystyle{\quad\mbox{if}\;\;N=3,4,5}, \\
\displaystyle{1 < p < \min \left(\frac{N+4-2\gamma}{N-2}, \frac{N+1}{N-3} \right)}&\displaystyle{\quad\mbox{if}\;\;N\geq6}\\.
 \end{array}
 \right.
$$
 Then, there exists $T>0$ depending only on the norm
$$ \|u_0\|_{H^1} +  \|u_1\|_{L^2} $$
and a unique solution $u$ to the problem $(\ref{1++})$ such that $u\in
C([0,T],H^1(\mathbb{R}^N))\cap C^1([0,T],L^2(\mathbb{R}^N)).$
\end{theorem}

Since
$$ \max\{p_1(N,\gamma),1/\gamma\} < 1 + \frac{4-2\gamma}{N-2}\quad\hbox{for all}\;\gamma\in(0,1),$$
the above local existence result enables one to extend the blow up result to all  values
of $p \in (0, \max\{p_1(N,\gamma),1/\gamma\} ]$ and all $\gamma \in (0,1).$

\begin{theorem}\label{Blow-up-theorem3in}
Let $N\geq3,$ $\gamma\in(0,1)$ and $1<p\leq  \max\{p_1(N,\gamma),1/\gamma\} .$ Assume that  $(u_0,u_1)\in H^1(\mathbb{R}^N)\times L^2(\mathbb{R}^N)$ is such that
$$
\int_{\mathbb{R}^N}u_0>0,\quad \int_{\mathbb{R}^N}u_1>0.
$$
 Then the solution of $(\ref{1++})$
blows up in finite time.\\
\end{theorem}

\noindent The organization of this paper is as follows. In Section
$\ref{sec2++},$ we give some properties, results and notations that will be used in the sequel.
In Section $\ref{sec3++},$ we present the local existence results of solutions for the
equation $(\ref{1++}).$ Section $\ref{sec4++},$ contains the blow-up results
of solutions to $(\ref{1++}).$


\section{Preliminaries, notations}\label{sec2++}
\setcounter{equation}{0}

In this section, we present some definitions, notations and results concerning the wave operator, fractional integrals and
fractional derivatives that will be used hereafter. For more information see \cite{GV3}, \cite{KTao}, \cite{KSTr} and \cite{SKM}.\\
Let us consider the inhomogeneous wave equation
\begin{equation}\label{inhomogeneous}
\left\{
\begin{array}{ll}
\displaystyle{u_{tt}-\Delta u=f,}&\displaystyle{(x,t)\in\mathbb{R}^N\times(0,T),}\\
\displaystyle{u(x,0)=u_0(x),\;u_t(x,0)=u_1(x),}&\displaystyle{x\in\mathbb{R}^N.}\\
\end{array}
\right.
\end{equation}
We define $K(t)$ and $\dot{K}(t)$ by $K(t):=\omega^{-1}\sin\omega t$ and $\dot{K}(t):=\cos\omega t$ where $\omega^{-1}$ is the inverse of the fractional laplacian operator $\omega:=(-\Delta)^{1/2}$ of order $1/2$ defined above. The solution of the Cauchy problem $(\ref{inhomogeneous})$ can be written, according to Duhamel's principle, as
\begin{equation}\label{eq.int1++}
   u(t) = \dot{K}(t)u_0+K(t)u_1 + \int^t_0 K(t-s)f(s)\,ds.
\end{equation}

The initial data $(u_0,u_1)$ of the problem $(\ref{inhomogeneous})$ will be taken in the energy space
\begin{equation}\label{Espaces}
\mathcal{H} = H^1(\mathbb{R}^N)\times L^2(\mathbb{R}^N)
\end{equation}
or more generally in
\begin{equation}\label{Espacesmu}
\mathcal{H}^\mu = H^\mu(\mathbb{R}^N)\times H^{\mu-1}(\mathbb{R}^N), \mu \geq 1.
\end{equation}

\noindent  We shall denote by $ \dot{H}^\mu(\mathbb{R}^N),$ $\mu \geq 0,$ the homogeneous Sobolev space of order $\mu \geq 0$ defined by
$$\begin{array}{ll}
\,\,\displaystyle{\dot{H}^\mu(\mathbb{R}^N)=\left\{u\in \mathcal{S}';\;(-\Delta)^{\mu/2}u\in L^2(\mathbb{R}^N)\right\}},
\end{array}
$$
where $\mathcal{S}'$ is the space of Schwartz' distributions and $(-\Delta)^{\mu/2}$ is the fractional laplacian operator defined by
$$(-\Delta)^{\mu/2}u(x):=\mathcal{F}^{-1}\left(|\xi|^\mu
\mathcal{F}(u)(\xi)\right)(x)$$ and $\mathcal{F}^{-1}$ stands the Fourier transform and its inverse, respectively.

The corresponding inhomogeneous Sobolev space $ H^\mu(\mathbb{R}^N)$ for any real $\mu$ is defined as
$$\begin{array}{ll}
\,\,\displaystyle{H^\mu(\mathbb{R}^N)=\left\{u\in \mathcal{S}';\;(1-\Delta)^{\mu/2}u\in L^2(\mathbb{R}^N)\right\}}.
\end{array}
$$

\noindent Next, we give the admissible version of the Strichartz estimates due to Keel and Tao \cite{KTao}. Before we state the theorem of Strichartz' estimates, we give the definition of $\sigma-$admissible pair where $\sigma=(N-1)/2$ for the wave equation.
\begin{definition}\label{definitionsigmaadmissible}$(\mbox{\cite[Definition~1.1]{KTao}})$
We say that the exponents pair $(q,r)$ is $\sigma-$admissible if $q,r\geq2,$ $(q,r,\sigma)\neq(2,\infty,1)$ and
\begin{equation}\label{sigmaadmissible}
\frac{1}{q}+\frac{\sigma}{r}\leq\frac{\sigma}{2}.
\end{equation}
If equality holds in $(\ref{sigmaadmissible}),$ we say that $(q,r)$ is sharp $\sigma-$admissible, otherwise we say that $(q,r)$ is nonsharp $\sigma-$ admissible. Note in particular that when $\sigma>1$ the endpoint
$$P=\left(2,\frac{2\sigma}{\sigma-1}\right)$$
is sharp $\sigma-$admissible. $\hfill\square$\\
\end{definition}
\begin{theorem}$(\mbox{\cite[Corollary~1.3]{KTao}})$
Suppose that $N\geq2$ and $(q,r)$ and $(\tilde{q},\tilde{r})$ are $(N-1)/2-$admissible pairs with $r,\tilde{r}<\infty.$ If $u$ is a $(\mbox{weak})$ solution to the problem $(\ref{inhomogeneous})$ in $\mathbb{R}^N\times[0,T]$ for some data $u_0\in H^\mu(\mathbb{R}^N), u_1\in H^{\mu-1}(\mathbb{R}^N),f\in L^{\tilde{q}'}([0,T];L^{\tilde{r}'}_x)$ and time $0<T<\infty,$ then
\begin{eqnarray}\label{Strichartzestimatehom}
\|u\|_{L^q([0,T];L^r_x)}&+&\|u\|_{C([0,T];\dot{H}^\mu)}+\|\partial_tu\|_{C([0,T];\dot{H}^{\mu-1})}\nonumber\\
&\leq& \overline{C}\left(\|u_0\|_{\dot{H}^\mu}+\|u_1\|_{\dot{H}^{\mu-1}}+\|f\|_{L^{\tilde{q}'}([0,T];L^{\tilde{r}'}_x)}\right),
\end{eqnarray}
under the assumption that the dimensional analysis $(\mbox{or "gap"})$ condition
\begin{equation}\label{gapcondition}
\frac{1}{q}+\frac{N}{r}=\frac{N}{2}-\mu=\frac{1}{\tilde{q}'}+\frac{N}{\tilde{r}'}-2
\end{equation}
holds, where $\overline{C}>0$ is a positive constant independent of $T.$ $\hfill\square$\\
\end{theorem}

\begin{rmk} In the above Theorem we denote by $\tilde{r}',\tilde{q}'$  the conjugate exponents of $\tilde{r},\tilde{q}$ and by $L^p_x:=L^p(\mathbb{R}^N)$ the standard Lebesgue $x$ space for all $1\leq p\leq\infty.$\\
The estimate $(\ref{Strichartzestimatehom})$ involves homogeneous Sobolev spaces. If we admit dependence of the constants on the length of the time interval $I=[T_1,T_2]$, taking the length $|I| = T_2-T_1 \leq 1$ and  $\mu \geq 1$ we can establish the inequality
\begin{eqnarray}\label{Strichartzestimate}
\|u\|_{L^q(I;L^r_x)}&+&\|u\|_{C(I;H^\mu)}+\|\partial_tu\|_{C(I;H^{\mu-1})}\nonumber\\
&\leq& C_0\left(\|u_0\|_{H^\mu}+\|u_1\|_{H^{\mu-1}}+\|f\|_{L^{\tilde{q}'}(I;L^{\tilde{r}'}_x)}\right),
\end{eqnarray}
where $C_0$ is independent of $|I| \leq 1.$ This inequality  is sufficient for the proof of local existence result and the existence of maximal interval of existence of the solution.

\end{rmk}

\begin{corollary}\label{Strichartzforu_0}$(\mbox{Strichartz estimates for $u_0$})$
Suppose that $N\geq2$ and $(q,r)$ is a $(N-1)/2-$admissible pair with $r<\infty.$ If $u_0\in H^\mu(\mathbb{R}^N)$, then
\begin{equation}\label{Strichartzestimateforu0}
\|\dot{K}(t)u_0\|_{L^q([0,T];L^r_x)}+\|\dot{K}(t)u_0\|_{C([0,T];H^1)}+\| \Delta K(t)u_0\|_{C([0,T];L^2)}\leq C \|u_0\|_{H^1},
\end{equation}
under the assumption that the  condition
\begin{equation}\label{gapconditionforu0}
\frac{1}{q}+\frac{N}{r}=\frac{N}{2}-1
\end{equation}
holds.  $\hfill\square$
\end{corollary}
\begin{corollary}\label{Strichartzforu_1}$(\mbox{Strichartz estimates for $u_1$})$
Suppose that $N\geq2$ and $(q,r)$ is a $(N-1)/2-$admissible pair with $r<\infty.$ If $u_1\in L^2(\mathbb{R}^N)$, then
\begin{equation}\label{Strichartzestimateforu1}
\|K(t)u_1\|_{L^q([0,T];L^r_x)}+\|K(t)u_1\|_{C([0,T];H^1)}+\|\dot{K}(t)u_1\|_{C([0,T];L^2)}\leq C \|u_1\|_{L^2},
\end{equation}
under the assumption that the gap condition
\begin{equation}\label{gapconditionforu1}
\frac{1}{q}+\frac{N}{r}=\frac{N}{2}-1
\end{equation}
holds. ${}\hfill\square$
\end{corollary}
\begin{corollary}\label{Strichartzforf}$(\mbox{Strichartz estimates for $f$})$
Suppose that $N\geq2$ and $(q,r)$ and $(\tilde{q},\tilde{r})$ are $(N-1)/2-$admissible pairs with $r,\tilde{r}<\infty.$ If $I=[T_1,T_2]$ is any time interval of length $|I|=T_2-T_1 \leq 1$ and  $f\in L^{\tilde{q}'}([0,T];L^{\tilde{r}'}_x),$  then
\begin{eqnarray}\label{Strichartzestimateforf}
\left\| \int_0^t K(t-s) f(s) ds\right\|_{L^q(I;L^r_x)}+\left\| \int_0^t K(t-s) f(s) ds\right\|_{C(I;H^1)}&+&\left\| \int_0^t \dot{K}(t-s) f(s) ds\right\|_{C(I;L^2)}\nonumber\\
&\leq& C_0 \|f\|_{L^{\tilde{q}'}(I;L^{\tilde{r}'}_x)},
\end{eqnarray}
under the assumption that the gap condition
\begin{equation}\label{gapconditionforf}
\frac{1}{q}+\frac{N}{r}=\frac{N}{2}-1=\frac{1}{\tilde{q}'}+\frac{N}{\tilde{r}'}-2
\end{equation}
holds.$\hfill\square$
\end{corollary}

Turning back to integral equation $(\ref{eq.int1++})$, we have to give a more precise definition of the
integral terms of the right hand side.

For the purpose we suppose that for some $T>0$ one can find admissible couple $(q,r)$  such that
the gap condition $(\ref{gapconditionforf})$ is satisfied and
$$ u \in X(T) = X_{q,r}(T) = C([0,T],H^1(\mathbb{R}^N))\cap C^1([0,T],L^2(\mathbb{R}^N)) \cap
L^{q}([0,T];L^{r}_x).$$

Then estimates of Corollary $\ref{Strichartzforu_0}$ and Corollary $\ref{Strichartzforu_1}$ guarantee that
$$ \dot{K}(t)u_0 \in X(T), K(t)u_1 \in X(T). $$

The estimate of Corollary $\ref{Strichartzforf}$ implies that
\begin{equation}\label{eq.XY}
   f \in Y(T) = Y_{\tilde{q}, \tilde{r}}(T) = L^{\tilde{q}'}([0,T];L^{\tilde{r}'}_x) \ \  \Longrightarrow \int^t_0 K(t-s)f(s)\,ds \in X(T)
\end{equation}
provided $(\tilde{q}, \tilde{r})$ is admissible and the gap condition $(\ref{gapconditionforf})$ is fulfilled.
Note that the integral in $(\ref{eq.XY})$ can be considered as Bochner integral in
$$ H^{-k}(\mathbb{R}^N) \supset  L^{\tilde{r}'}_x $$ due to the Sobolev embedding with
$$ \frac{1}{\tilde{r}'} - \frac{1}{2} = \frac{k}{N}.$$

\noindent The final part of this section is devoted to some basic properties of Riemann-Liouville
fractional derivatives. If $AC[0,T]$ is the space of all functions which are
absolutely continuous on $[0,T]$ with $0<T<\infty,$ then, for $f\in
AC[0,T],$ the left-handed and right-handed Riemann-Liouville
fractional derivatives $D^{\alpha}_{0|t}f(t)$ and
$D^{\alpha}_{t|T}f(t)$ of order $\alpha\in(0,1)$ are defined by (see
\cite{KSTr})
\begin{eqnarray}
D^\alpha_{0|t}f(t)&:=&DJ^{1-\alpha}_{0|t}f(t)\label{estiF1},\\
D^\alpha_{t|T}f(t)&:=&-\frac{1}{\Gamma(1-\alpha)}D\int_t^T(s-t)^{-\alpha}f(s)\,ds,
\end{eqnarray}
for all $t\in[0,T],$ where $D:=\frac{d}{dt}$ and
\begin{equation}\label{estiF2++}
J^\alpha_{0|t}g(t):=\frac{1}{\Gamma(\alpha)}\int_0^t(t-s)^{\alpha-1}g(s)\,ds
\end{equation}
is the Riemann-Liouville fractional integral (see \cite{KSTr}), for all $g\in L^q(0,T)$ $(1\leq q\leq\infty).$\\
Furthermore, for every $f,g\in C([0,T]),$ such that
$D^\alpha_{0|t}f(t),D^\alpha_{t|T}g(t)$ exist and are continuous,
for all $t\in[0,T],$ $0<\alpha<1,$ we have the formula of
integration by parts (see  $(2.64)$ p. 46 in \cite{SKM})
\begin{equation}\label{F.D++}
\int_0^T \left(D^\alpha_{0|t}f\right)(t)g(t)\,dt \;=\; \int_0^T
f(t)\left(D^\alpha_{t|T}g\right)(t)\,dt.
\end{equation}
Note also that, for all $f\in AC^{n+1}[0,T]$ and all integer $n\geq0,$ we have (see $(2.2.30)$ in
\cite{KSTr})
\begin{equation}\label{B.0++}
(-1)^{n}D^n.D^\alpha_{t|T}f=D^{n+\alpha}_{t|T}f,
\end{equation}
where
$$AC^{n+1}[0,T]:=\left\{f:[0,T]\rightarrow\mathbb{R}\;\hbox{and}\;D^nf\in
AC[0,T]\right\}$$
and $D^n$ is the usual $n$ times derivative.\\
 Moreover, for all $1\leq q\leq\infty,$ the
following formula (see \cite[Lemma~2.4 p.74]{KSTr})
\begin{equation}\label{estiF3++}
    D^\alpha_{0|t}J^\alpha_{0|t}=Id_{L^q(0,T)}
\end{equation}
holds almost everywhere on $[0,T].$\\
Later on, we will use the following results.\\
\noindent If $w_1(t)=\left(1-{t}/{T}\right)_+^\sigma,$
$t\geq0,$ $T>0,$ $\sigma\gg1,$ then
\begin{eqnarray}
D_{t|T}^\alpha
w_1(t)&=&\frac{(1-\alpha+\sigma)\Gamma(\sigma+1)}{
\Gamma(2-\alpha+\sigma)}T^{-\sigma}(T-t)_+^{\sigma-\alpha},\label{2.14++}\\
D_{t|T}^{\alpha+1}
w_1(t)&=&\frac{(1-\alpha+\sigma)(\sigma-\alpha)
\Gamma(\sigma+1)}{\Gamma(2-\alpha+\sigma)}T^{-\sigma}(T-t)_+^{\sigma-\alpha-1},\label{110++}\\
D_{t|T}^{\alpha+2}
w_1(t)&=&\frac{(1-\alpha+\sigma)(\sigma-\alpha)(\sigma-\alpha-1)\Gamma(\sigma+1)}{\Gamma(2-\alpha+\sigma)}T^{-\sigma}(T-t)_+^{\sigma-\alpha-2},\label{1000++}
\end{eqnarray}
for all $\alpha\in(0,1);$ so
\begin{equation}\label{estiF4++}
\left(D_{t|T}^\alpha w_1\right)(T)=0\hspace{4 mm};\hspace{4
mm}\left(D_{t|T}^\alpha w_1\right)(0)=C\;T^{-\alpha},
\end{equation}
and
\begin{equation}\label{estiF40++}
\left(D_{t|T}^{\alpha+1} w_1\right)(T)=0\hspace{4 mm};\hspace{4
mm}\left(D_{t|T}^{\alpha+1} w_1\right)(0)=\tilde{C}\;T^{-\alpha-1},
\end{equation}
where
$$C=\frac{(1-\alpha+\sigma)\Gamma(\sigma+1)}{\Gamma(2-\alpha+\sigma)}\quad\mbox{and}\quad \tilde{C}=\frac{(1-\alpha+\sigma)(\sigma-\alpha)\Gamma(\sigma+1)}{\Gamma(2-\alpha+\sigma)}.$$
Indeed, using the Euler change of variable $y={(s-t)}/{(T-t)},$ we
get
\begin{eqnarray*}
D_{t|T}^\alpha w_1(t) &:=& -\frac{1}{\Gamma(1-\alpha)}D
\left[\int_t^T(s-t)^{-\alpha}\left(1-\frac{s}{T}\right)^\sigma\,ds\right] \\
&=& -\frac{T^{-\sigma}}{\Gamma(1-\alpha)}D\left[(T-t)^{1-\alpha+\sigma}
\int_0^1(y)^{-\alpha}(1-y)^\sigma\,ds\right] \\
&=&+\frac{(1-\alpha+\sigma)B(1-\alpha;\sigma+1)}{\Gamma(1-\alpha)}
T^{-\sigma}(T-t)^{\sigma-\alpha},
        \end{eqnarray*}
where $B(\cdotp;\cdotp)$ stands for the beta function. Then, $(\ref{2.14++})$
follows using the relation
$$B(1-\alpha;\sigma+1)=\frac{\Gamma(1-\alpha)\Gamma(\sigma+1)}{\Gamma(2-\alpha+\sigma)}.$$
Furthermore, $(\ref{110++})$ and $(\ref{1000++})$ follow from the formula $(\ref{B.0++})$ applied to $(\ref{2.14++}).$ $\hfill\square$\\


\section{Local existence and uniqueness   theorems for mild and weak solutions}\label{sec3++}
\setcounter{equation}{0}

 First we recall the definition of local
mild solution for the problem $(\ref{1++}).$
\begin{definition}$(\mbox{Mild solution of $(\ref{1++})$})$
 Given any $\mu\geq1$ and any $T>0$ we say that $$u \in C([0,T],H^\mu(\mathbb{R}^N))\cap C^1([0,T],H^{\mu-1}(\mathbb{R}^N)) $$ is a mild solution of $(\ref{1++})$ with initial data
$$ (u_0,u_1) \in \mathcal{H}^\mu $$
 if $u$ satisfies the integral equation
\begin{equation}\label{mildsolution1++}
u(t)=\dot{K}(t)u_0+K(t)u_1+ \int_0^t K(t-s) J^\alpha_{0|s}(|u|^{p}))(s) ds,\quad t\in[0,T].
\end{equation}
\end{definition}

\begin{definition} \label{def:weak}$(\mbox{Weak solution of $(\ref{1++})$})$
\ Given any $T>0$ we say that $u$ is a weak solution of $(\ref{1++})$ if there exist admissible couples
$(q,r)$ and $(\tilde{q}, \tilde{r})$ so that the gap condition $(\ref{gapconditionforf})$ is fulfilled,
$$u\in
X(T) = C([0,T],H^1(\mathbb{R}^N))\cap C^1([0,T],L^2(\mathbb{R}^N)) \cap L^{q}([0,T];L^{r}_x),$$
$$J^\alpha_{0|t}(|u|^{p}))(t) \in Y(T) = Y_{\tilde{q}, \tilde{r}}(T) = L^{\tilde{q}'}([0,T];L^{\tilde{r}'}_x)$$
 and $u$ satisfies the integral equation
\begin{equation}\label{wesolution1++}
u(t)=\dot{K}(t)u_0+K(t)u_1+ \int_0^t K(t-s) J^\alpha_{0|s}(|u|^{p}))(s) ds,\quad t\in[0,T].
\end{equation}
\end{definition}

Our first goal of this section is to establish the existence and uniqueness of mild solutions.

\begin{theorem}\label{Theoremmild0++}$(\mbox{local existence of unique mild solution of $(\ref{1++})$})$\\
Suppose $(u_0,u_1)\in H^1(\mathbb{R}^N) \times L^2(\mathbb{R}^N),$ $N\geq 1,$  $\gamma\in(0,1)$ and let $p>1$ be such that
\begin{equation}\label{conditionoverp}
\left\{
\begin{array}{ll}
\displaystyle{1 < p \leq \frac{N}{N-2}}&\displaystyle{\quad\mbox{if}\;\;N>2}\\
\displaystyle{1 < p < \infty }&\displaystyle{\quad\mbox{if}\;\;N=1,2}. \\
 \end{array}
 \right.
\end{equation}
 Then, there exist $T>0$ depending only on the norm
$$ \|u_0\|_{H^1} +  \|u_1\|_{L^2} $$
and a unique mild solution $u$ to the problem $(\ref{1++})$ such that $u\in
C([0,T],H^1(\mathbb{R}^N))\cap C^1([0,T],L^2(\mathbb{R}^N)).$
\end{theorem}
\begin{pf} For any $N \geq 1$ we apply the energy estimate
$$\|u\|_{C([0,T];H^{1})}+\|\partial_tu\|_{C([0,T];L^2)} \leq C_0 \left(\|u_0\|_{\dot{H}^1}+\|u_1\|_{L^2}\right) + C_0 \|J^\alpha_{0|t}(|u|^{p})(t) \|_{L^{1}([0,T];L^2_x)} . $$
Here and below $C_0=C_0(T)$ remains bounded, when $0 \leq T \leq 1.$
Then we have to show the estimates
\begin{equation}\label{eq.con1mild}
   \|J^\alpha_{0|t}(|u|^{p}))(t) \|_{L^{1}([0,T];L^2_x)}\leq C(T)\|u\|^p_{C([0,T];H^1)}
\end{equation}
and
\begin{equation}\label{eq.con2mild}
   \|J^\alpha_{0|t}(|u|^{p}))(t)- J^\alpha_{0|t}(|v|^{p}))(t) \|_{L^{1}([0,T];L^2_x)}\leq C(T)\|u-v\|_{C([0,T];H^1)}
   \left( \|u\|^{p-1}_{C([0,T];H^1)} + \|v\|^{p-1}_{C([0,T];H^1)}\right)
\end{equation}
with some constant $C(T)$ satisfying the property
$$ \lim_{T \rightarrow 0} C(T) = 0.$$
Once these estimates are established an application of a contraction mapping principle in $$X(T)
= C([0,T],H^1(\mathbb{R}^N))\cap C^1([0,T],L^2(\mathbb{R}^N))$$ will complete the proof.

We shall verify only $(\ref{eq.con1mild})$, since the proof of  $(\ref{eq.con2mild})$ is similar. We have
$$\|J^\alpha_{0|t}(|u|^{p})\|_{L^{1}([0,T];L^{2}_x)}\leq \|J^\alpha_{0|t}(\|u\|_{L^{2p}_x}^{p})\|_{L^{1}([0,T])}.$$
 For $N=1,2$ we have the Sobolev embedding
 \begin{equation}\label{sobolevimbedding0}
H^1(\mathbb{R}^N)\hookrightarrow L^{2p}(\mathbb{R}^N),
\end{equation}
valid for $2 < 2p < \infty.$ For $N \geq 3$ we have the same embedding provided
  the condition $p \leq N/(N-2)$ is fulfilled. Hence,
we get
\begin{equation}\label{eq2++mild}
\|J^\alpha_{0|t}(|u|^{p})\|_{L^{1}([0,T];L^{2}_x)} \leq C_1^p\|J^\alpha_{0|t}(\|u\|^p_{H^1(\mathbb{R}^N)})\|_{L^{1}([0,T])},
\end{equation}
where $C_1$ is the positive constant of the Sobolev imbedding. Using the fact that $u\in X(T),$ we have
$$\|J^\alpha_{0|t}(\|u\|^p_{H^1(\mathbb{R}^N)})\|_{L^{1}([0,T])}\leq \frac{1}{(2-\gamma)\Gamma(2-\gamma)}T^{2-\gamma}\|u\|^p_{C([0,T];H^1)}.$$
This completes the check of $(\ref{eq.con1mild})$ and  the proof of the Theorem.
 $\hfill\square$\\
\end{pf}

To get local mild solution for some $p > N/(N-2)$ we have to impose different assumptions on $N,p.$

\begin{theorem}\label{TheoremmildP0++}$(\mbox{local existence of unique mild solution of $(\ref{1++})$})$\\
Suppose $N\geq3,$  $\gamma\in(0,1)$ and let $p>N/(N-2)$ be such that
\begin{equation}\label{conditionoverp00}
\begin{array}{ll}
\displaystyle{p^2 - \frac{pN}{2}+ \frac{N}{2} \geq  0 }. \\
 \end{array}
\end{equation}
If
$(u_0,u_1)\in H^\mu(\mathbb{R}^N) \times H^{\mu - 1}(\mathbb{R}^N),$ where
$\mu = N/2 - 1/(p-1) > 1.$
 Then, there exists $T>0$ depending only on the norm
$$ \|u_0\|_{H^\mu} +  \|u_1\|_{H^{\mu-1}} $$
and a unique mild solution $u$ to the problem $(\ref{1++})$ such that $u\in
C([0,T],H^\mu(\mathbb{R}^N))\cap C^1([0,T],H^{\mu-1}(\mathbb{R}^N)).$
\end{theorem}
\begin{pf} We follow the proof of the previous result and take
$$ q=\infty, \ r= \frac{2N}{(N-2\mu)}.$$ Using the Sobolev embedding with some $\mu >1,$ we get
\begin{equation}\label{eq.sob12}
    \|u\|_{L^{\infty}([0,T];L^{r}_x)} \leq C \|(-\Delta)^{(\mu-1)/2}u\|_{L^{\infty}([0,T];L^{r_1}_x)} \leq
C\|u\|_{C([0,T];\dot{H}^{\mu})},
\end{equation}
where
$$ r_1= \frac{2N}{N-2}.$$ These Sobolev embeddings are fulfilled because
\begin{equation}\label{eq.par}
    \frac{1}{r_1} - \frac{1}{r} = \frac{\mu-1}{N}, \ \ \frac{1}{2} - \frac{1}{r_1} = \frac{1}{N}
\end{equation}
and in the second inequality in $(\ref{eq.sob12})$ we use the classical Sobolev inequality
$$ \|f\|_{L^{r_1}_x} \leq C \|f\|_{\dot{H}^{1}}$$
with $f= (-\Delta)^{(\mu-1)/2}u.$

\noindent Note that $ v=(-\Delta)^{(\mu-1)/2}u$ is a solution to the equation
$$  v_{tt}-\Delta
v=D_t^{\gamma - 1}(-\Delta_x)^{(\mu -1)/2}\left(|u(t)|^p\right),$$
so applying the classical energy estimate for this wave equation we find
$$\|u\|_{C([0,T];\dot{H}^{\mu})}+\|\partial_tu\|_{C([0,T];H^{\mu-1})} \leq C_0 \left(\|u_0\|_{\dot{H}^\mu}+\|u_1\|_{\dot{H}^{\mu-1}}\right) + C_0 \|J^\alpha_{0|t}((-\Delta)^{(\mu-1)/2}|u|^{p}))(t) \|_{L^{1}([0,T];L^2_x)} . $$
From $(\ref{eq.sob12})$ we conclude
\begin{eqnarray} \label{enermu}
 \|u\|_{L^{\infty}([0,T];L^{r}_x)} + \|(-\Delta)^{(\mu-1)/2}u\|_{L^{\infty}([0,T];L^{r_1}_x)} &+&\|u\|_{C([0,T];\dot{H}^{\mu})}+  \\ \nonumber
 \|\partial_tu\|_{C([0,T];H^{\mu-1})} \leq C \left(\|u_0\|_{\dot{H}^\mu}+\|u_1\|_{\dot{H}^{\mu-1}}\right) &+& C \|J^\alpha_{0|t}((-\Delta)^{(\mu-1)/2}|u|^{p}))(t) \|_{L^{1}([0,T];L^2_x)} .
\end{eqnarray}

Now we are in position to apply the following inequality (see for example Lemma 2.3 in \cite{GV2} or \cite{Run}, \cite{RuSi})
$$ \| |u|^p \|_{H^{\mu - 1}} \leq C \|(-\Delta)^{(\mu-1)/2}u\|_{L^{r_1}} \| u \|^{p-1}_{L^{r_2(p-1)}},$$
where $1/r_1+1/r_2= 1/2$ and $n/r_1 >\mu -1.$ Note that our choice of $r_1$ implies $r_2=N$ so we can use the relation
$$ r= \frac{2N}{(N-2\mu)} , \ \mu = \frac{N}{2}- \frac{1}{p-1} \ \Longrightarrow \ N(p-1) = r.$$
It is important to notice that the above estimate of the nonlinear term $|u|^p$ is valid only for $\mu -1 \leq p,$
since $p>1$ might be not integer. The inequality $\mu -1 \leq p,$ as well our choice of $\mu$ lead to the inequality
$$ p^2 - \frac{pN}{2}+ \frac{N}{2} \geq  0.$$

We can proceed further as in the proof of the previous Theorem and we can show the estimates
\begin{equation}\label{eq.con1mildP}
   \|J^\alpha_{0|t}((-\Delta)^{(\mu-1)/2}|u|^{p}))(t) \|_{L^{1}([0,T];L^2_x)}\leq C(T)\|u\|^p_{C([0,T];H^\mu)}
\end{equation}
and
\begin{equation}\label{eq.con2mildP}
   \|J^\alpha_{0|t}((-\Delta)^{(\mu-1)/2}|u|^{p}))(t)- J^\alpha_{0|t}((-\Delta)^{(\mu-1)/2}|v|^{p}))(t) \|_{L^{1}([0,T];L^2_x)}\leq C(T)\|u-v\|_{C([0,T];H^\mu)}
   \left( \|u\|^{p-1}_{C([0,T];H^\mu)} + \|v\|^{p-1}_{C([0,T];H^\mu)}\right)
\end{equation}
where  $C(T)$ is an increasing, continuous in  $(0,1]$ function, satisfying the property
$$ \lim_{T \rightarrow 0} C(T) = 0.$$
Once these estimates are established an application of a contraction mapping principle in $$X(T)
= C([0,T],H^\mu(\mathbb{R}^N))\cap C^1([0,T],H^{\mu - 1}(\mathbb{R}^N))$$ and this completes the proof. $\hfill\square$\\
\end{pf}

\begin{rmk}
The condition
$$ p^2 - \frac{pN}{2}+ \frac{N}{2} \geq  0$$
is automatically satisfied if $3 \leq N \leq 8.$ The condition becomes very restrictive in the case of space dimensions
$9 \leq N \leq 20.$ One can show that the critical exponent $p_2(N)$ is strictly smaller than any of the roots of
 $p^2 - \frac{pN}{2}+ \frac{N}{2} =  0$ is $N$ is large enough, namely $N \geq 20.$

\end{rmk}

\begin{rmk}
 The proof of the existence and uniqueness of mild solutions is done without Strichartz' estimate, using only the energy estimate
 $$\|u\|_{C([0,T];\dot{H}^1)}+\|\partial_tu\|_{C([0,T];L^2)} \leq \overline{C} \left(\|u_0\|_{\dot{H}^\mu}+\|u_1\|_{\dot{H}^{\mu-1}}+\|f\|_{L^{1}([0,T];L^{2}_x)}\right)  $$
 and Sobolev embedding. For this the restrictive assumption of type
 $$p^2 - \frac{pN}{2}+ \frac{N}{2} \geq  0$$
 can not be avoided. Nevertheless, one can prove the existence of a maximal
time $0<T_{\max}\leq\infty$ and a unique mild solution $u$ to the problem $(\ref{1++})$ such that $u\in C([0,T_{\max});H^\mu(\mathbb{R}^N))\cap  C^1([0,T_{\max}); H^{\mu-1}(\mathbb{R}^N)).$ Moreover, if $T_{\max}<\infty,$ we have
$$(\|u(t)\|_{H^\mu(\mathbb{R}^N)}+\|u_t(t)\|_{H^{\mu-1}(\mathbb{R}^N)})\longrightarrow\infty\quad \mbox{as}\;\;t\rightarrow T_{\max}.
$$
Furthermore, if
\begin{equation}\label{condition2}
\mbox{supp}u_i\subset B(r):=\{x\in\mathbb{R}^n:\;|x|< r\},\quad r>0,\;i=0,1,
\end{equation}
$u(t,\cdotp)$ is supported in the ball $B(t+r).$
We note that, we can extend our local existence theorem to the case $N\geq1$ by assuming that the initial data satisfies furthermore $(\ref{condition2})$ and using the fact that $A$ is a skew-adjoint operator in $H^1\times L^2$ (see \cite[Theorem~6.2.2, p.~76]{CH}) instead to use Strichartz' estimate.
\end{rmk}

\noindent To cover completely the case $N/(N-2) < p \leq (N+2)/(N-2)$ and show that the problem $(\ref{1++})$ is locally well posed in $H^1$
one has to use effectively the Strichartz estimate ( as it is done in \cite{GV1}, \cite{GV2} ) and work with weak solutions
of Definition $\ref{def:weak}.$ In this work we need local existence and existence of maximal time interval for the solution,
while in in \cite{GV1}, \cite{GV2}   the global Cauchy problem is studied. For this we can prove that the problem $(\ref{1++})$ is locally well posed for a larger interval $ p \in (1,  \min \{{(N+4-2\gamma)}/{(N-2)}, {(N+1)}/{(N-3)} \}).$

\begin{theorem}\label{Theoremweak0++}$(\mbox{local existence of unique weak solution of $(\ref{1++})$})$\\
Given $(u_0,u_1)\in H^1(\mathbb{R}^N) \times L^2(\mathbb{R}^N),$ $N\geq 1,$  $\gamma\in(0,1)$ and let $p>1$ be such that
\begin{equation}\label{conditionoverp0}
\left\{
\begin{array}{lll}
\displaystyle{1 < p < \infty }&\displaystyle{\quad\mbox{if}\;\;N=1,2}, \\
\displaystyle{1 < p < \frac{N+4-2\gamma}{N-2}}&\displaystyle{\quad\mbox{if}\;\;N=3,4,5}, \\
\displaystyle{1 < p < \min \left(\frac{N+4-2\gamma}{N-2}, \frac{N+1}{N-3} \right)}&\displaystyle{\quad\mbox{if}\;\;N\geq6}\\.
 \end{array}
 \right.
\end{equation}
 Then, there exist $T>0$ depending only on the norm
$$ \|u_0\|_{H^1} +  \|u_1\|_{L^2} $$
and a unique weak solution $u$ to the problem $(\ref{1++})$ such that $u\in
C([0,T],H^1(\mathbb{R}^N))\cap C^1([0,T],L^2(\mathbb{R}^N)).$
\end{theorem}
\begin{pf} We shall consider only the case $N > 2,$ since for $N=1,2$ we already have established the existence of mild solutions.
There is no lack of generality if we suppose
$$ \frac{N}{N-2} < p < \frac{N+4-2\gamma}{N-2},$$
 since for
$$ 1 < p \leq \frac{N}{N-2}$$ Theorem $\ref{Theoremmild0++}$ guarantees that local mild solution exists and it is unique.

We take the following admissible couple
\begin{equation}\label{eq.strchoice}
    \frac{1}{r} = \min \left(\frac{N+1}{2pN}, \frac{N-2}{2N} -\varepsilon \right), \ \frac{1}{q}= \frac{N-2}{2} - \frac{N}{r}
\end{equation}
with $\varepsilon > 0$ small enough.

To explain how we arrived at this choice and then how to complete the proof of the Theorem, we write the general conditions of admissibility as well as the gap condition
\begin{equation}\label{eq.admis11}
    \frac{1}{q} + \frac{N-1}{2r} \leq \frac{N-1}{4} , \ \ \frac{1}{\tilde{q}} + \frac{N-1}{2\tilde{r}} \leq \frac{N-1}{4},
\end{equation}
\begin{equation}\label{gap11}
   \frac{1}{q} = \frac{N-2}{2} - \frac{N}{r} , \ \ \frac{1}{\tilde{q}^\prime} = \frac{N+2}{2} - \frac{N}{\tilde{r}^\prime}.
\end{equation}

To apply a contraction mapping principle we need to apply Strichartz estimate as well as the estimate
\begin{equation}\label{eq.con1a}
   \|D_t^{-\alpha}(|u|^{p}))(t) \|_{L^{\tilde{q}^\prime}([0,T];L^{\tilde{r}^\prime}_x)}\leq C(T)\|u\|^p_{L^q([0,T];L^r_x)}.
\end{equation}
For this we take
\begin{equation}\label{eqball}
    p \tilde{r}^\prime = r.
\end{equation}

The Sobolev embedding $\dot{H}_{q^*}^{\alpha}(0,T) \subset L^{\tilde{q}^\prime}(0,T)$ with
$$ \frac{1}{q^*} = \frac{1}{\tilde{q}^\prime} + \alpha$$
combined with the H\"older inequality imply
$$
\|D_t^{-\alpha}(|u|^{p}))(t) \|_{L^{\tilde{q}^\prime}([0,T]}\leq C \||u|^{p}(t) \|_{L^{q^*}(0,T)} \leq C(T)\|u\|^p_{L^q(0,T)}
$$
with $\lim_{T \rightarrow 0} C(T)=0$ provided
$ q^* p < q$ i.e.
\begin{equation}\label{eq.mainds}
    \frac{p}{q} < \frac{1}{\tilde{q}^\prime} + \alpha.
\end{equation}
If we take the gap condition $(\ref{gap11})$ and the relation $(\ref{eqball})$ we see that
we are able to express the parameters $q, \tilde{q}^\prime $ and $ \tilde{r}^\prime$ as functions of $r,p.$
$$ \frac{1}{q} = \frac{N-2}{2} - \frac{N}{r},\ \frac{1}{\tilde{q}^\prime} = \frac{N+2}{2} - \frac{pN}{r}, \ \frac{1}{\tilde{r}^\prime} = \frac{p}{r}. $$
Substitution in $(\ref{eq.mainds})$ leads to the inequality
$$  p < \frac{N+4-2\gamma}{N-2}$$ while admissible  conditions and natural requirements
$ 1 < \tilde{q^\prime} \leq 2  \leq q < \infty$ can be rewritten as
$$ \frac{1}{2p} < \frac{1}{r} \leq \frac{N+1}{2pN}, \ \frac{N-3}{2N} \leq \frac{1}{r} < \frac{N-2}{2N}.$$
This domain is non empty if and only if
$$ \frac{1}{2p} < \frac{N-2}{2N}, \ \frac{N-3}{2N} < \frac{N+1}{2pN}.$$
Since we already made the assumption $p > N/(N-2)$ we see that
$$ p < \frac{N+1}{N-3}
$$ has to be imposed too.

This observation suggests the choice $(\ref{eq.strchoice})$ with $\varepsilon >0$ so small that the domain
 $$ \frac{1}{2p} < \frac{1}{r} \leq \frac{N+1}{2pN}, \ \frac{N-3}{2N} \leq \frac{1}{r} < \frac{N-2}{2N}$$
is nonempty and it is sufficient to apply contraction principle.

 \noindent To be more precise we have to prove the existence and uniqueness of the fixed point for the integral equation
$$ u(t) = \dot{K}(t)u_0 + K(t)u_1 + \int_0^t K(t-s) D_t^{-\alpha}|u|^p(s) ds $$
such that
$$ u(t) \in X(T) = C([0,T],H^1(\mathbb{R}^N))\cap C^1([0,T],L^2(\mathbb{R}^N)) \cap L^{q}([0,T];L^{r}_x).$$

Applying the Strichartz estimate $(\ref{Strichartzestimate})$
 as in the proof of Theorem $\ref{Theoremmild0++},$ we
 obtain estimate
  \begin{equation}\label{eq.timeex}
   \left\|\dot{K}(t)u_0 + K(t)u_1 + \int_0^t K(t-s) D_t^{-\alpha}|u|^p(s) ds \right\|_{X(T)} \leq C_0 \|(u_0,u_1)\|_{\mathcal{H}} + C(T) \|u\|^p_{X(T)}
  \end{equation}
 and
 $$ \left\|\int_0^t K(t-s) \left( D_t^{-\alpha}|u|^p(s) - D_t^{-\alpha}|v|^p(s) \right) ds \right\|_{X(T)} \leq
 C(T) \|u-v\|_{X(T)} \left( \|u\|^{p-1}_{X(T)} +  \|v\|^{p-1}_{X(T)}\right) ,  $$
 where
 $ \lim_{T \rightarrow 0} C(T) = 0.$

 Applying the contraction principle we get existence and uniqueness of weak solution. The fact that the time interval depends only on the energy norm $$  \|(u_0,u_1)\|_{\mathcal{H}}$$ of the initial data follows directly from $(\ref{eq.timeex})$ since
 the fixed point $u \in X(T)$ will satisfies the estimate
 $$ \left\|u \right\|_{X(T)} \leq C_0 \|(u_0,u_1)\|_{\mathcal{H}} + C(T) \|u\|^p_{X(T)} $$
 and this estimate implies
 $$ \left\|u \right\|_{X(T)} \leq 2 C_0 \|(u_0,u_1)\|_{\mathcal{H}} $$
 if
 $$ C(T) 2^p \left(C_0 \|(u_0,u_1)\|_{\mathcal{H}} \right)^{p-1} < 1.$$

 This complete the proof of the theorem. $\hfill\square$\\

\end{pf}

\begin{rmk}
 Since the time interval $[0,T]$ depends only on the size of the energy norm of the initial data, one can prove the existence of a maximal
time $0<T_{\max}\leq\infty$ and a unique weak solution $u$ to the problem $(\ref{1++})$ such that $u\in C([0,T_{\max});H^1(\mathbb{R}^N))\cap  C^1([0,T_{\max}); L^2(\mathbb{R}^N)).$ Moreover, if $T_{\max}<\infty,$ we have
$$(\|u(t)\|_{H^1(\mathbb{R}^N)}+\|u_t(t)\|_{L^2(\mathbb{R}^N)})\longrightarrow\infty\quad \mbox{as}\;\;t \nearrow T_{\max}.
$$
Furthermore, if
\begin{equation}\label{condition2w}
\mbox{supp}u_i\subset B(r):=\{x\in\mathbb{R}^n:\;|x|< r\},\quad r>0,\;i=0,1,
\end{equation}
$u(t,\cdotp)$ is supported in the ball $B(t+r).$

\end{rmk}

Since
$$ \min \left(\frac{N+4-2\gamma}{N-2}, \frac{N+1}{N-3} \right) = \frac{N+1}{N-3}  $$
for $0 < \gamma < (N-5)/(N-3)$ we have the following.

\begin{corollary}
\label{cor:loc} Suppose $N \geq 6$,
$$ 0 < \gamma < \frac{N-5}{N-3}$$
 and $$ 1 < p < \frac{N+1}{N-3}.$$
 Then, for any $(u_0,u_1) \in H^1(\mathbb{R}^N) \times L^2(\mathbb{R}^N)$ there exists a maximal time $0<T_{\max}\leq\infty$ and a unique solution $u$ to the problem $(\ref{1++})$ such that $u\in
C([0,T_{\max}),H^1(\mathbb{R}^N))\cap C^1([0,T_{\max}),L^2(\mathbb{R}^N)).$
\end{corollary}


\section{Blow-up theorems}\label{sec4++}
\setcounter{equation}{0}

This section is devoted to the blow-up of solutions of the problem $(\ref{1++}),$ assuming initial data are in the energy space
 $$ (u_0,u_1) \in \mathcal{H}$$ and $p,$ $ \gamma$ satisfy appropriate subcritical inequalities.
 To do this, we have to introduce the definition of the solution of $(\ref{1++})$ in distributional sense and to prove that the mild and weak solutions of $(\ref{1++})$  are  solutions in distributional sense of the same equation, because our blow up argument  is based on this fact.\\

\begin{rmk}
As we shall use solutions  in distributional sense, the natural question is why we discussed  mild and weak solutions? The answer is the following property of weak and mild solutions: either $T_{\max}=\infty$ or else $T_{\max}<\infty$ and $\|u(t)\|_{H^1(\mathbb{R}^N)}+\|u_t(t)\|_{L^2(\mathbb{R}^N)}\rightarrow\infty$ as $t\rightarrow T_{\max}.$
\end{rmk}
\begin{definition}$(\mbox{Solution in distributional sense})$
Let $u_0,u_1\in L_{Loc}^1(\mathbb{R}^N).$ We say that $u$ is a  solution of $(\ref{1++})$ in distributional sense, if and only if $u\in L^p((0,T),L_{Loc}^p(\mathbb{R}^N))$ satisfies
\begin{eqnarray}\label{weaksolution++}
&{}&\int_0^T\int_{\Omega}J^{\alpha}_{0|t}(|u|^{p})(x,t)\varphi(x,t)+\int_{\Omega}u_1(x)\varphi(x,0)-\int_{\Omega}u_0(x)\varphi_t(x,0)\nonumber\\
&{}&=\int_0^T\int_{\Omega}u(x,t)\varphi_{tt}(x,t)-\int_0^T\int_{\Omega}u(x,t)\Delta\varphi(x,t)
\end{eqnarray}
for all compactly supported function $\varphi\in C^2([0,T]\times\mathbb{R}^N)$ such that $\varphi(\cdotp,T)=0$ and $\varphi_t(\cdotp,T)=0,$ where $\alpha:=1-\gamma\in(0,1),$ $\Omega:=$supp$\varphi.$
\end{definition}
\begin{lemma}\label{MildWeak++}$(\mbox{Mild \ or \ Weak $\rightarrow$ Distributional})$
Assume that $(u_0,u_1)\in \mathcal{H}$ and $\gamma\in(0,1).$ Let $u$ be the mild or weak solution of $(\ref{1++})$ and let $p>1$ satisfies  $(\ref{conditionoverp})$ or $(\ref{conditionoverp0})$ respectively with $(u,u_t)\in C([0,T],\mathcal{H}),$ then $u$ is a distributional solution of $(\ref{1++}),$ for all $T>0.$
\end{lemma}
\begin{pf} We shall consider the case of mild solutions, since the argument works as well for the weak solutions. Let $T>0,$ $u$ be a mild solution of $(\ref{1++})$ and $\varphi\in C^2([0,T]\times\mathbb{R}^N)$ be a compactly supported function such that $\varphi(\cdotp,T)=0,$ $\varphi_t(\cdotp,T)=0$ and supp$\varphi=:\Omega.$ Then, $u$ is a fixed point
for the integral equation
\begin{equation}\label{mildB++}
u(t)= \dot{K}(t)u_0+\Delta K(t)u_1+(K\ast J^\alpha_{0|t}(|u|^{p}))(t),
\end{equation}
and we have
\begin{equation}\label{mildA++}
u_t(t)=\Delta K(t)u_0+\dot{K}(t)u_1+(\dot{K}\ast J^\alpha_{0|t}(|u|^{p}))(t).
\end{equation}
So, after multiplying $(\ref{mildA++})$ by $\varphi$ and integrating over $\mathbb{R}^N,$ we obtain
\begin{eqnarray*}
 \int_{\Omega}u_t(x,t)\varphi(x,t)&=& \int_{\Omega}\Delta K(t)u_0(x)\varphi(x,t)+ \int_{\Omega}\dot{K}(t)u_1(x)\varphi(x,t)\\
 &+&\int_{\Omega}\int_0^t\dot{K}(t-s) J^\alpha_{0|s}(|u|^{p})(x,s)\,ds\varphi(x,t).
    \end{eqnarray*}
Then
 \begin{eqnarray}\label{differentiation++}
\frac{d}{dt} \int_{\Omega}u_t(x,t)\varphi(x,t)&=&\frac{d}{dt} \int_{\Omega}\Delta K(t)u_0(x)\varphi(x,t)+\frac{d}{dt} \int_{\Omega}\dot{K}(t)u_1(x)\varphi(x,t)\nonumber\\
 &+&\int_{\Omega}\frac{d}{dt}\int_0^t\dot{K}(t-s) J^\alpha_{0|s}(|u|^{p})(x,s)\,ds\varphi(x,t).
    \end{eqnarray}
Now, using the fact that the Laplacian is a negative self-adjoint operator,  we have:
\begin{eqnarray}\label{int1++}
&{}& \frac{d}{dt} \int_{\Omega}\Delta K(t)u_0(x)\varphi(x,t)+\frac{d}{dt} \int_{\Omega}\dot{K}(t)u_1(x)\varphi(x,t)\nonumber\\
&{}&=\int_{\Omega} \Delta \left[\dot{K}(t)u_0(x)+ K(t)u_1(x)\right]\varphi(x,t)+ \int_{\Omega}\left[\Delta K(t)u_0(x)+\dot{K}(t)u_1(x)\right]\varphi_t(x,t)\nonumber\\
&{}&=\int_{\Omega}\left[\dot{K}(t)u_0(x)+ K(t)u_1(x)\right]\Delta\varphi(x,t)+ \int_{\Omega}\left[\Delta K(t)u_0(x)+\dot{K}(t)u_1(x)\right]\varphi_t(x,t),\nonumber\\
&{}&
\end{eqnarray}
and
\begin{eqnarray}\label{int2++}
&{}&\int_{\Omega}\frac{d}{dt}\int_0^t\dot{K}(t-s) f(x,s)\,ds\varphi(x,t)\nonumber\\
&{}&=\int_{\Omega}f(x,t)\varphi(x,t)+\int_{\Omega}\int_0^t \Delta\left(K(t-s)f(x,s)\right)\,ds\varphi(x,t)\nonumber\\
&{}&+ \int_{\Omega} \int_0^t\dot{K}(t-s)f(x,s)\,ds\varphi_t(x,t)\nonumber\\
 &{}&=  \int_{\Omega}f(x,t)\varphi(x,t)+\int_{\mathbb{R}^N}(K\ast f)(x,t)\Delta\varphi(x,t)+ \int_{\Omega} (\dot{K}\ast f)(x,t)\varphi_t(x,t)\quad
     \end{eqnarray}
where $f:=J_{0|t}^\alpha
    \left(|u|^{p}\right)\in C([0,T];L^2(\Omega)).$\\
    Thus, using $(\ref{mildB++})-(\ref{mildA++})$ and  $(\ref{int1++})-(\ref{int2++}),$ we conclude that $(\ref{differentiation++})$ implies that
\begin{eqnarray}\label{final1++}
 \frac{d}{dt} \int_{\Omega}u_t(x,t)\varphi(x,t)&=& \int_{\Omega} u(x,t)\Delta\varphi(x,t)+\int_{\Omega} u_t(x,t)\varphi_t(x,t)\nonumber\\
 &+&\int_{\Omega}f(x,t)\varphi(x,t).
    \end{eqnarray}
Next, after integrating in time $(\ref{final1++})$ over $[0,T]$ and using the fact that $\varphi(\cdotp,T)=0$ and $\varphi_t(\cdotp,T)=0,$ we conclude that
\begin{eqnarray}\label{final2++}
 -\int_{\Omega}u_1(x)\varphi(x,0)&=&\int_0^T \int_{\Omega} u(x,t)\Delta\varphi(x,t)-\int_0^T\int_{\Omega} u(x,t)\varphi_{tt}(x,t)\nonumber\\
 &-&\int_{\Omega}u_0(x)\varphi_t(x,0)+\int_{\Omega}f(x,t)\varphi(x,t).
    \end{eqnarray}
$\hfill\square$
\end{pf}
As
$$p_1=p_2=\frac{1}{\gamma}={N}/{(N-2)}\Longleftrightarrow \gamma={(N-2)}/{N},$$
so we have to distinguish two cases: $\gamma>{(N-2)}/{N}$ and $\gamma\leq{(N-2)}/{N}.$ Moreover, in the case when $\gamma>{(N-2)}/{N},$ we note that ${1}/{\gamma}<p_2<p_1<{N}/{(N-2)}$ for $N=2m,$ $m\in\mathbb{N}\setminus\{0,1\},$ and ${1}/{\gamma}<p_1<p_2<{N}/{(N-2)}$ for $N=2m+1,$ $m\in\mathbb{N}^*,$ while ${N}{(N-2)}<p_2<p_1<{1}/{\gamma}<$ when $\gamma<{(N-2)}/{N}.$  For that, we have the following blow-up theorems.
\begin{theorem}\label{blowuptheorem++}$(\gamma>{(N-2)}/{N}\;\mbox{and $N\in\{1\}\cup\{2m,\;m\in\mathbb{N}^*\}$})$\\
Let $1<p\leq {N}/{(N-2)},$ if $N\geq 3,$ and $p\in(1,\infty),$ if $N=1,2.$ Assume that $N\in\{1\}\cup\{2m,\;m\in\mathbb{N}^*\},$  ${(N-2)}/{N}<\gamma<1$ and $(u_0,u_1)\in H^1(\mathbb{R}^N)\times L^2(\mathbb{R}^N)$ such that
\begin{equation}\label{averages}
\int_{\mathbb{R}^N}u_0>0,\quad \int_{\mathbb{R}^N}u_1>0.
\end{equation}
If $p\leq p_1,$ where $p_1$ is given in $(\ref{conditionoverp1}),$ then the solution of $(\ref{1++})$
blows up in finite time.\\
\end{theorem}
\begin{pf} The proof proceeds by contradiction. Let $u$ be a global
mild solution of the problem $(\ref{1++}),$ then $u$ is a  mild solution of $(\ref{1++})$ in $C([0,T],H^1(\mathbb{R}^N))\cap C^1([0,T],L^2(\mathbb{R}^N))$ for all $T\gg1.$ Using Lemma $\ref{MildWeak++},$ we have
\begin{eqnarray}\label{weaksolution0++sup}
&{}&\int_0^T\int_{\mbox{supp$\varphi$}}J^{\alpha}_{0|t}(|u|^p)(x,t)\varphi(x,t)+\int_{\mbox{supp$\varphi$}}u_1(x)\varphi(x,0)\,dx-\int_{\mbox{supp$\varphi$}}u_0(x)\varphi_t(x,0)\nonumber\\
&{}&=\int_0^T\int_{\mbox{supp$\varphi$}}u(x,t)\varphi_{tt}(x,t)-\int_0^T\int_{\mbox{supp$\varphi$}}u(x,t)\Delta\varphi(x,t)
\end{eqnarray}
for all compactly supported function $\varphi\in C^2([0,T]\times\mathbb{R}^N)$ such that $\varphi(\cdotp,T)=0$ and $\varphi_t(\cdotp,T)=0,$ where $\alpha:=1-\gamma\in(0,1).$\\
Now, we have to distinguish two cases:\\

\noindent $\bullet$ \underline{The case $p<p_1$}: Let $\varphi(x,t)=D^\alpha_{t|T}\left(\tilde{\varphi}(x,t)\right):=D^\alpha_{t|T}\left(\left(\varphi_1(x)\right)^\ell\varphi_2(t)\right)$ with
$\varphi_1(x):=\Phi\left({|x|}/{T}\right),$
$\varphi_2(t):=\left(1-{t}/{T}\right)^\eta_+,$ where $\ell,\eta\gg1$ and $\Phi$ be a
smooth non-increasing function such that
$$\Phi(r)=\left\{\begin {array}{ll}\displaystyle{1}&\displaystyle{\quad\mbox{if }0\leq r\leq 1,}\\
\displaystyle{0}&\displaystyle{\quad\mbox {if }r\geq 2,}
\end {array}\right.$$\\
$0\leq \Phi \leq 1,$ $|\Phi^{'}(r)|\leq C_1/r,$ for all $r>0.$ Then, we have
\begin{eqnarray}\label{weaksolution0++}
&{}&\int_{\Omega_T}J^{\alpha}_{0|t}(|u|^p)(x,t)D^\alpha_{t|T}\tilde{\varphi}(x,t)+\int_{\Omega}u_1(x)D^\alpha_{t|T}\tilde{\varphi}(x,0)-\int_{\Omega}u_0(x)DD^\alpha_{t|T}\tilde{\varphi}(x,0)\nonumber\\
&{}&=\int_{\Omega_T}u(x,t)D^2D^\alpha_{t|T}\tilde{\varphi}(x,t)-\int_{\Omega_T}u(x,t)\Delta D^\alpha_{t|T}\tilde{\varphi}(x,t),
\end{eqnarray}
where
$$\Omega_T=[0,T]\times \Omega\;\;\mbox{for}\;\;\Omega:=\left\{x\in{\mathbb{R}}^N\hspace{2
mm};\hspace{2mm}|x|\leq 2 T\right\},\;\;\int_{\Omega_T}=\int_{\Omega_T}\,dx\,dt,\;\;\int_{\Omega}=\int_{\Omega}\,dx.$$
Moreover, from $(\ref{F.D++}),(\ref{B.0++}),(\ref{estiF4++})$ and $(\ref{estiF40++})$ we
may write
\begin{eqnarray}\label{estiAC++}
&{}&\int_{\Omega_T}D^{\alpha}_{0|t}J^{\alpha}_{0|t}(|u|^p)\;\tilde{\varphi}
\;+\;C\;T^{-\alpha}\int_{\Omega}\left(\varphi_1(x)\right)^\ell u_1(x)\;+\;\tilde{C}\;T^{-\alpha-1}\int_{\Omega}\left(\varphi_1(x)\right)^\ell u_0(x)\nonumber\\
&{}&=\;
\int_{\Omega_T}u\left(\varphi_1(x)\right)^\ell D^{2+\alpha}_{t|T}\varphi_2(t)\;+\;\int_{\Omega_T}u(-\Delta_x)\left(\varphi_1(x)\right)^{\ell}D^\alpha_{t|T}\varphi_2(t).
\end{eqnarray}
So, $(\ref{estiF3++})$ and the formula $\Delta\left(\varphi_1^\ell\right)=\ell\varphi_1^{\ell-1}\Delta\varphi_1+\ell(\ell-1)\varphi_1^{\ell-2}|\nabla\varphi_1|^2$ will allow us to write:
\begin{eqnarray}\label{new2.15++}
&{}&\int_{\Omega_T}|u|^p\;\tilde{\varphi}\;+\;C\;T^{-\alpha}\int_{\Omega}\left(\varphi_1(x)\right)^\ell u_1(x)\;+\;\tilde{C}\;T^{-\alpha-1}\int_{\Omega}\left(\varphi_1(x)\right)^\ell u_0(x)\nonumber\\
&{}&=\;\int_{\Omega_T}u\left(\varphi_1(x)\right)^\ell D^{2+\alpha}_{t|T}\varphi_2(t)\;- \;C\int_{\Omega_T}u\left(\varphi_1(x)\right)^{\ell-1}
\Delta_x\varphi_1(x)\;D^\alpha_{t|T}\varphi_2(t)\nonumber\\
&{}&\quad\;-\;C\int_{\Omega_T}u\left(\varphi_1(x)\right)^{\ell-2}
|\nabla\varphi_1(x)|^2\;D^\alpha_{t|T}\varphi_2(t)\nonumber\\
&{}&\leq\;
\int_{\Omega_T}|u|\left(\varphi_1(x)\right)^\ell \left|D^{2+\alpha}_{t|T}\varphi_2(t)\right|\;+\;C\int_{\Omega_T}|u|\left(\varphi_1(x)\right)^{\ell-1}
\left|\Delta_x\varphi_1(x)\;D^\alpha_{t|T}\varphi_2(t)\right|\nonumber\\
&{}&\quad+\;C\int_{\Omega_T}|u|\left(\varphi_1(x)\right)^{\ell-2}
|\nabla\varphi_1(x)|^2\;\left|D^\alpha_{t|T}\varphi_2(t)\right|
\end{eqnarray}
Therefore, as the condition $(\ref{averages})$ implies
\begin{equation}\label{positivity1++}
\int_{\mbox{supp$\varphi_1$}}\left(\varphi_1(x)\right)^\ell u_0(x)\geq0,\quad\int_{\mbox{supp$\varphi_1$}}\left(\varphi_1(x)\right)^\ell u_1(x)\geq0,
\end{equation}
(here supp$\varphi_1=\Omega),$ we obtain
\begin{eqnarray}\label{2.15++}
\int_{\Omega_T}|u|^p\;\tilde{\varphi}&\leq& \int_{\Omega_T}|u|\;\tilde{\varphi}^{1/p}\tilde{\varphi}^{-1/p}\left(\varphi_1(x)\right)^\ell
\left|D^{2+\alpha}_{t|T}\varphi_2(t)\right|\nonumber\\
&{}&\;+\;C\int_{\Omega_T}|u|\;\tilde{\varphi}^{1/p}\tilde{\varphi}^{-1/p}\;\left(\varphi_1(x)\right)^{\ell-1}
\left|\Delta_x\varphi_1(x)\;D^\alpha_{t|T}\varphi_2(t)\right|\nonumber\\
&{}&\;+\;C\int_{\Omega_T}|u|\;\tilde{\varphi}^{1/p}\tilde{\varphi}^{-1/p}\;\left(\varphi_1(x)\right)^{\ell-2}
|\nabla\varphi_1(x)|^2\;\left|D^\alpha_{t|T}\varphi_2(t)\right|,
\end{eqnarray}
So, using the Young inequality
\begin{equation}\label{2.16++}
ab\leq\;\frac{1}{3p}a^{\;p}\;+\;\frac{3^{\tilde{p}-1}}{\tilde{p}}\;
b^{\;\tilde{p}}\qquad\hbox{where}\;\; p\tilde{p}=p+\tilde{p},\;\;p>1,\tilde{p}>1,\;\;a>0,b>0,
\end{equation}
with
$$\left\{\begin{array}{l}
a=|u|\;\tilde{\varphi}^{1/p},\\
b=\tilde{\varphi}^{-1/p}\left(\varphi_1(x)\right)^\ell
\left|D^{2+\alpha}_{t|T}\varphi_2(t)\right|,\\
\end{array}
\right.
$$
in the first integral of the right hand side of $(\ref{2.15++}),$
$$\left\{\begin{array}{l}
a=|u|\;\tilde{\varphi}^{1/p},\\
b=C\;\tilde{\varphi}^{-1/p}\;\left(\varphi_1(x)\right)^{\ell-1}
\left|\Delta_x\varphi_1(x)\;D^\alpha_{t|T}\varphi_2(t)\right|,\\
\end{array}
\right.
$$
in the second integral of the right hand side of $(\ref{2.15++})$ and with
$$\left\{\begin{array}{l}
a=|u|\;\tilde{\varphi}^{1/p},\\
b=C\;\tilde{\varphi}^{-1/p}\;\left(\varphi_1(x)\right)^{\ell-2}
|\nabla\varphi_1(x)|^2\left|D^\alpha_{t|T}\varphi_2(t)\right|,\\
\end{array}
\right.
$$
in the third integral of the right hand side of $(\ref{2.15++}),$ we
obtain
\begin{eqnarray}\label{esti2++}
 \int_{\Omega_T}|u(x,t)|^p\;\tilde{\varphi}(x,t)&\leq& C\int_{\Omega_T}\left(\varphi_1\right)^\ell
\left(\varphi_2\right)^{-\frac{1}{p-1}}
\left|D^{2+\alpha}_{t|T}\varphi_2\right|^{\tilde{p}}\nonumber\\
&{}&+\;C\int_{\Omega_T}\left(\varphi_1\right)^{\ell-\tilde{p}}
\left(\varphi_2\right)^{-\frac{1}{p-1}}\left|\Delta_x\varphi_1D^\alpha_{t|T}\varphi_2\right|^{\tilde{p}}\nonumber\\
&{}&+\;C\int_{\Omega_T}\left(\varphi_1\right)^{\ell-2\tilde{p}}
\left(\varphi_2\right)^{-\frac{1}{p-1}}|\nabla\varphi_1|^{2\tilde{p}}\left|D^\alpha_{t|T}\varphi_2\right|^{\tilde{p}}.
\end{eqnarray}
\noindent At this stage, we introduce the scaled variables: $\tau=
T^{-1}t$ and $\xi= T^{-1}x;$ using formulas $(\ref{2.14++})$ and $(\ref{1000++})$ in the right
hand-side of $(\ref{esti2++}),$ we obtain:
\begin{equation}\label{2.17++}
\int_{\Omega_T}|u(x,t)|^p\;\tilde{\varphi}(x,t) \leq
\;C\;T^{-\delta},
\end{equation}
\noindent where $\delta:=(2+\alpha)
\tilde{p}-1-N,$ $C=C(|\Omega_1 |\;,\; |\Omega_2
|), (|\Omega_i |$ stands for the measure of $\Omega_i,$ for $
i=1,2),$ with
$$\Omega_1:=\left\{\xi\in{\mathbb{R}}^N\hspace{2
mm};\hspace{2mm}|\xi|\leq 2
\right\}\quad,\quad\Omega_2:=\left\{\tau\geq0\hspace{2
mm};\hspace{2mm}\tau\leq 1\right\}.$$\\
Passing to the limit in $(\ref{2.17++}),$  as $T$ goes to $\infty,$ and taking into account the fact that $p<p_1$ $(\Longleftrightarrow \delta>0),$ we
conclude that
$$\lim_{T\rightarrow\infty}\int_0^{T}\int_{|x|\leq
2T}|u(x,t)|^p\;\tilde{\varphi}(x,t)\,dx\,dt =0.$$
Using the dominated convergence theorem,  we infer that
$$\int_0^\infty\int_{\mathbb{R}^N}|u(x,t)|^p\,dx\,dt =0\qquad\;\Longrightarrow\qquad\; u=0\;\;\mbox{for all $t$ and a.e. $x$}.$$
This contradicts the fact that $\int_{\mathbb{R}^N}u_0>0.$\\

\noindent $\bullet$ \underline{The case $p=p_1$}: In this case, we take $\tilde{\varphi}(x,t)=\left(\varphi_1(x)\right)^\ell\varphi_2(t)$ with
$\varphi_1(x):=\Phi\left({|x|}/{B^{-1}T}\right),$
$\varphi_2(t):=\left(1-{t}/{T}\right)^\eta_+,$ instead of the one used in the last case, where $\ell,\eta\gg1$ and $1\leq B<T$ large enough such that when $T\rightarrow\infty,$ we don't have $B\rightarrow\infty$ at the same time. Here $\Phi$ is the same function used above.\\
So, by repeating the same computations as in the case $p<p_1,$ we obtain
\begin{eqnarray}\label{1001++}
&{}&\int_{\Sigma_B}|u|^p\;\tilde{\varphi}
\;+\;C\;T^{-\alpha}\int_{\Omega_B}\left(\varphi_1(x)\right)^\ell u_1(x)\;+\;\tilde{C}\;T^{-\alpha-1}\int_{\Omega_B}\left(\varphi_1(x)\right)^\ell u_0(x)\nonumber\\
&{}&\leq\;\int_{\Sigma_B}|u|\;\tilde{\varphi}^{1/p}\tilde{\varphi}^{-1/p}\left(\varphi_1(x)\right)^\ell
\left|D^{2+\alpha}_{t|T}\varphi_2(t)\right|\nonumber\\
&{}&\quad+\;C\int_{\Delta_B}|u|\;\tilde{\varphi}^{1/p}\tilde{\varphi}^{-1/p}\;\left(\varphi_1(x)\right)^{\ell-1}
\left|\Delta_x\varphi_1(x)\;D^\alpha_{t|T}\varphi_2(t)\right|\nonumber\\
&{}&\quad+\;C\int_{\Delta_B}|u|\;\tilde{\varphi}^{1/p}\tilde{\varphi}^{-1/p}\;\left(\varphi_1(x)\right)^{\ell-2}
|\nabla\varphi_1(x)|^2\left|D^\alpha_{t|T}\varphi_2(t)\right|,
\end{eqnarray}
where
$$\Sigma_B=[0,T]\times \Omega_{B}:=[0,T]\times\left\{x\in{\mathbb{R}}^N\hspace{2
mm};\hspace{2mm}|x|\leq 2B^{-1} T\right\},\;\;\int_{\Sigma_B}=\int_{\Sigma_B}\,dx\,dt,\;\;\int_{\Omega_B}=\int_{\Omega_B}\,dx$$
and
$$\Delta_B:=[0,T]\times \left\{x\in{\mathbb{R}}^N\hspace{2
mm};\hspace{2mm}B^{-1}T \leq |x|\leq 2
B^{-1}T\right\},\;\;\int_{\Delta_B}=\int_{\Delta_B}\,dx\,dt.$$
Moreover, using the Young inequality
\begin{equation}\label{1002++}
ab\leq\;\frac{1}{p}a^{\;p}\;+\;\frac{1}{\tilde{p}}\;b^{\;\tilde{p}}\qquad\hbox{where}\;\; p\tilde{p}=p+\tilde{p},\;\;p>1,\tilde{p}>1,\;\;a>0,b>0,
\end{equation}
with
$$\left\{\begin{array}{l}
a=|u|\;\tilde{\varphi}^{1/p},\\
b=\tilde{\varphi}^{-1/p}\left(\varphi_1(x)\right)^\ell
\left|D^{2+\alpha}_{t|T}\varphi_2(t)\right|,\\
\end{array}
\right.
$$
in the first integral of the right hand side of $(\ref{1001++}),$ and
using H\"older's inequality
$$\int_{\Delta_B}ab\leq\;\left(\int_{\Delta_B}a^{\;p}\right)^{1/p}\;
\left(\int_{\Delta_B}b^{\;\tilde{p}}\right)^{1/\tilde{p}},\;\;p>1,\tilde{p}>1,\;\;a>0,b>0,$$
with
$$\left\{\begin{array}{l}
a=|u|\;\tilde{\varphi}^{1/p},\\
b=\tilde{\varphi}^{-1/p}\;\left(\varphi_1(x)\right)^{\ell-1}
\left|\Delta_x\varphi_1(x)\;D^\alpha_{t|T}\varphi_2(t)\right|,\\
\end{array}
\right.
$$
in the second integral of the right hand side of $(\ref{1001++})$ and with
$$\left\{\begin{array}{l}
a=|u|\;\tilde{\varphi}^{1/p},\\
b=\tilde{\varphi}^{-1/p}\;\left(\varphi_1(x)\right)^{\ell-2}
|\nabla\varphi_1(x)|^2\left|D^\alpha_{t|T}\varphi_2(t)\right|,\\
\end{array}
\right.
$$
in the third integral of the right hand side of $(\ref{1001++}),$ and taking account of $(\ref{positivity1++}),$ we obtain
\begin{eqnarray}\label{1003++}
&{}&\int_{\Sigma_B}|u(x,t)|^p\;\tilde{\varphi}(x,t)\nonumber\\
&{}&\leq C\int_{\Sigma_B}\left(\varphi_1\right)^\ell
\left(\varphi_2\right)^{-\frac{1}{p-1}}
\left|D^{2+\alpha}_{t|T}\varphi_2\right|^{\tilde{p}}\nonumber\\
&{}&\quad+\;C\left(\int_{\Delta_B}|u|^p\;\tilde{\varphi}\right)^{1/p}\left(\int_{\Delta_B}
\left(\varphi_1\right)^{\ell-\tilde{p}}
\left(\varphi_2\right)^{-\frac{1}{p-1}}\left|\Delta_x
\varphi_1D^\alpha_{t|T}\varphi_2\right|^{\tilde{p}}\right)^{1/\tilde{p}}\nonumber\\
&{}&\quad+\;C\left(\int_{\Delta_B}|u|^p\;\tilde{\varphi}\right)^{1/p}\left(\int_{\Delta_B}
\left(\varphi_1\right)^{\ell-2\tilde{p}}
\left(\varphi_2\right)^{-\frac{1}{p-1}}|\nabla\varphi_1|^{2\tilde{p}}\left|D^\alpha_{t|T}\varphi_2\right|^{\tilde{p}}\right)^{1/\tilde{p}}.
\end{eqnarray}
\noindent Taking account of the scaled variables: $\tau= T^{-1}t,$
 $\xi= \left( T/B \right)^{-1}x,$ the formulas $(\ref{2.14++}),(\ref{1000++})$ and the fact that $p=p_1,$ we get
\begin{equation}\label{104++}
\int_{\Sigma_B}|u(x,t)|^p\;\tilde{\varphi}(x,t)\leq
\;C\;B^{-N}\;+\;C\;B^{2-\frac{N}{\tilde{p}}}\left(\int_{\Delta_B}|u(x,t)|^p\;\tilde{\varphi}(x,t)\right)^{1/p}.
\end{equation}
Now, from $(\ref{2.17++})$ and the fact that $(p=p_1\Longleftrightarrow \delta=0),$ we have the following implication
$$\lim_{T\rightarrow\infty}\int_{\Sigma_B}|u(x,t)|^p\;\tilde{\varphi}(x,t)\leq
\;C\;\;\Longrightarrow\;\;\int_0^\infty\int_{\mathbb{R}^N}|u(x,t)|^p\leq
\;C,$$
and so
$$
\lim_{T\rightarrow\infty}\left(\int_{\Delta_B}|u|^p\;\tilde{\varphi}\right)^{1/p}=\left(\lim_{T\rightarrow\infty}\int_0^T\int_{|x|\leq 2B^{-1}T}|u|^p\;\tilde{\varphi}-\lim_{T\rightarrow\infty}\int_0^T\int_{|x|\leq B^{-1}T}|u|^p\;\tilde{\varphi} \right)^{1/p}=0.$$
Thus, passing to the limit in $(\ref{104++}),$ as
$T\rightarrow\infty,$ we get
$$\int_0^\infty\int_{\mathbb{R^N}}|u(x,t)|^p\,dx\,dt \leq
\;C\;B^{-N}.$$ Then, taking the limit when $B$ goes to
infinity, we obtain $u=0$ for all $t$ and for almost every $x;$ contradiction with the fact that $\int_{\mathbb{R}^N}u_0>0.$ $\hfill\square$\\
\end{pf}

\begin{theorem}\label{Blow-up-theorem2}$(\gamma>{(N-2)}/{N}\;\mbox{and $N=2m+1,$ $m\in\mathbb{N}^*$})$\\
Let $1<p\leq {N}/{(N-2)},$ $N=2m+1,$ $m\in\mathbb{N}^*,$ ${(N-2)}/{N}<\gamma<1$ for $N=3$ and $\max\{1-(p-1)(N-3)/2,{(N-2)}/{N}\}<\gamma<1$ for $N>3.$ Assume that $(u_0,u_1)\in H^1(\mathbb{R}^N)\times L^2(\mathbb{R}^N)$ satisfy $(\ref{condition2})$ such that
\begin{equation}\label{averages1}
 \int_{\mathbb{R}^N}u_1>0,\quad\mbox{for $N=3$}\qquad\mbox{and}\qquad \int_{\mathbb{R}^N}|x|^{-1}u_0>0,\; \int_{\mathbb{R}^N}u_1>0,\quad\mbox{for $N>3$}.
\end{equation}
If $p< p_2,$ where $p_2$ is given in $(\ref{p_2}),$ then the solution of $(\ref{1++})$
blows up in finite time.\\
\end{theorem}
\begin{pf} The first step is to obtain a differential inequality. Let $u$ be the mild solution of the problem $(\ref{1++}).$ Using the proof of Lemma $\ref{MildWeak++},$ we have
\begin{equation}\label{weaksolution0}
\frac{d^2}{dt}\int_{\mbox{supp$\varphi$}}u(x,t)\varphi(x,t)\,dx=\int_{\mbox{supp$\varphi$}}u(x,t)\Delta\varphi(x,t)\,dx+\int_{\mbox{supp$\varphi$}}J^{\alpha}_{0|t}(|u|^p)(x,t)\varphi(x,t) \,dx
\end{equation}
for all $0\leq t<T_{\max}$ and all compactly supported function $\varphi\in C^2(\mathbb{R}^N).$  Fix $0<T_0<T_{\max}$ and take $\varphi\in C^2(\mathbb{R}^N)$ with $\varphi\equiv1$ on $B(r+T_0).$ Then, for all $0\leq t\leq T_0,$ $(\ref{weaksolution0})$ implies
\begin{equation}\label{weaksolution1}
\frac{d^2}{dt}\int_{\mathbb{R}^N}u(x,t)\,dx=\int_{\mathbb{R}^N}J^{\alpha}_{0|t}(|u|^p)(x,t)\,dx.
\end{equation}
Actually, equation $(\ref{weaksolution1})$ holds on $[0,T_{\max})$ since $T_0$ was arbitrary.\\
Now, due to the positivity of the operator $K$ only in three dimension, we have to study two cases.\\
\noindent $\bullet$ \underline{The case $N=3$}: For $r\leq t<T_{\max}$ (if $T_{\max}\leq r$ there is nothing to prove) define
\begin{equation}\label{firstdefinitionofF}
F(t)=\int_{\mathbb{R}^3}u(x,t)\,dx.
\end{equation}
Using the compact support of $u(\cdotp,t)$ and H\"older's inequality, it follows from $(\ref{weaksolution1})$ and $(\ref{firstdefinitionofF})$ that
\begin{equation}\label{weaksolution2}
\ddot{F}(t)\geq J^{\alpha}_{0|t}[(r+\cdotp)^{-3(p-1)}|F(\cdotp)|^p](t).
\end{equation}
For details, see \cite{John}. On the other hand, it is well known that the operator $K$ in the integral equation $(\ref{mildsolution1++})$ is positive. Therefore, $(\ref{mildsolution1++})$ implies that
 \begin{equation}\label{solution3}
u(x,t)\geq v(x,t),
\end{equation}
where $v:= \dot{K}(t)u_0+K(t)u_1$. Since $({d^2}/{dt^2})\int_{\mathbb{R}^3}v(x,t)\,dx=0,$ we have
 \begin{equation}\label{solution4}
\int_{\mathbb{R}^3}v(x,t)\,dx=C_{u_1}t+C_{u_0},
\end{equation}
where $C_{u_i}:=\int u_i\,dx,$ $i=0,1.$ Using the strong Huygen's principle, we have
\begin{equation}\label{solution5}
\mbox{supp}v(x,t)\subset\{t-r<|x|<t+r\},\quad t>r.
\end{equation}
Combining $(\ref{solution3})-(\ref{solution5}),$ using H\"older's inequality, one has
\begin{equation}\label{solution6}
C_{u_1}t+C_{u_0}\leq C(t+r)^{2(p-1)/p}\left(\int_{\mathbb{R}^3}|u(x,t)|^p\,dx\right)^{1/p}.
\end{equation}
Next, as in $(\ref{weaksolution2}),$ we obtain from $(\ref{solution6})$
 \begin{equation}\label{weaksolution3}
\ddot{F}(t)\geq J^{\alpha}_{0|t}\left[ \int_{\mathbb{R}^3}|u(x,\cdotp)|^p\,dx\right](t)\geq J^{\alpha}_{0|t}[(C_{u_1}t+C_{u_0})^p(r+\cdotp)^{-2(p-1)}](t)\geq J^{\alpha}_{0|t}(Ct^{-(p-2)})=Ct^{\alpha-(p-2)},
\end{equation}
where we have used the condition $(\ref{averages1})$, for $t$ large. Integrating twice, one has
 \begin{equation}\label{weaksolution4}
F(t)\geq Ct^{2+\alpha-(p-2)}\geq (r+t)^{\alpha_1},\qquad t\;\mbox{large},
\end{equation}
where $\alpha_1:=2+\alpha-(p-2).$ Turning back to $(\ref{weaksolution2})$ we can get
after integration twice
$$ F(t) \geq  C (r+t)^{\alpha_2}, $$
where
$$ \alpha_2 = p \alpha_1 - 3(p-1) + \alpha + 2.$$
Generally, we can write
$$ F(t) \geq  C (r+t)^{\alpha_k}, $$
where
$$ \alpha_{k+1} = p \alpha_k - 3(p-1) + \alpha + 2.$$
To assure that this sequence is increasing we need
$$\alpha_2 > \alpha_1$$ and a simple calculation shows that this is equivalent to $(\ref{p_2}).$ This is exactly the condition that means that $p>1$ is subcritical i.e. $p < p_2.$
Once the condition $ \alpha_2 > \alpha_1$ is verified one can verify that $\alpha_k$ tends to $\infty$ and deduce
the following estimates
\begin{equation}\label{eq.unN}
  F(t) \geq C_N (t+r)^N , \ \ \forall N \geq 1.
\end{equation}
\begin{equation}\label{eq.unNFsec}
  \ddot{F}(t) \geq C \int_0^t (t-s)^{\alpha - 1} F(s)^{p_1} ds , \ \ 1 < p_1 < p.
\end{equation}
Now we are in position to apply appropriate modification of  \cite[Lemma 4]{Sideris} and conclude that $T_{max} < \infty.$
\begin{lemma}
If $F(t) \in C^2([0,T))$ is an increasing positive function that satisfies  $(\ref{eq.unN})$ and $(\ref{eq.unNFsec})$ with
some $p_1 > 1 > \alpha > 0.$ Then $T < \infty.$
\end{lemma}
\begin{pf} Set
$$ G(t) = \int_0^t (t-s)^{\beta} F(s) ds,$$
where $\beta$ is a positive number such that
$$ \beta > \frac{\alpha}{(p_1-1)} - 1.$$
Then applying the H\"older inequality, we find
$$ G(t) \leq \left( \int_0^t (t-s)^{\alpha+\beta} F(s)^{p_1} ds\right)^{1/p_1} \left(\int_0^t (t-s)^{\beta  -\alpha/(p_1-1)} ds \right)^{(p_1-1)/p_1} \leq $$ $$\leq C (t+r)^{(\beta +1)(p_1-1)/p_1 -\alpha/p_1} \left( \int_0^t (t-s)^{\alpha-1} F(s)^{p_1} ds \right)^{1/p_1}, $$
since $ \beta  -\alpha/(p_1-1) > -1.$ Hence,
$$ \int_0^t (t-s)^{\alpha + \beta} F(s)^{p_1} ds  \geq C (t+r)^{\alpha - (\beta+1)( p_1-1)} \  G(t)^{p_1}, $$
and applying the estimate
$$ \ddot{G}(t) = \int_0^t (t-\tau)^{\beta} \ddot{F}(\tau) d\tau \geq
\int_0^t \int_0^\tau (t-\tau)^\beta (\tau - s)^{\alpha -1} F(s)^{p_1} ds d\tau \ \geq $$
$$ \geq C \int_0^t (t-s)^{\alpha + \beta} F(s)^{p_1} ds \geq C (t+r)^{\alpha - (\beta+1)( p_1-1)} \  G(t)^{p_1}.$$
The estimates
$$ \ddot{G}(t)  \geq C (t+r)^{\alpha - (\beta+1)( p_1-1)} \  G(t)^{p_1}$$
and
$$ G(t) \geq C_N (t+r)^N , \ \ \forall N \geq 1 $$
enables one to apply \cite[Lemma 4]{Sideris} and conclude that $T < \infty.$
This completes the proof of the Lemma.
\end{pf}

\noindent $\bullet$ \underline{The case $N>3$}: Let
$$F(t)=\int_0^t(t-s)^{(N-5)/2}\int_{\mathbb{R}^N}u(s,x)\,dx\,ds,\qquad r\leq t<T_{\max}.$$
We know that in the case $N=3$ the kernel $K$ is positive while in the high dimension space $N>3$ is not. So, we follow the approach of Sideris \cite{Sideris} and  defined  $F(\cdotp)$ with the purpose to use \cite[Lemma 5]{Sideris} and get the positivity.\\
Differentiating $F(t)$ twice and using $(\ref{weaksolution1}),$ we obtain
$$\ddot{F}(t)=\frac{N-5}{2}t^{(N-7)/2}C_{u_0}+t^{(N-5)/2}C_{u_1}+\int_0^t(t-s)^{(N-5)/2}\int_0^s(s-\sigma)^{\alpha-1}
\int_{\mathbb{R}^N}|u(\sigma,x)|^p\,dx\,d\sigma\,ds.$$
For $t$ large, inverting the order of integration and then using the compact support of $u(\cdotp,t),$ we get
\begin{eqnarray}\label{estimationA}
\ddot{F}(t)&\geq& \int_0^t(t-s)^{(N-5)/2}\int_0^s(s-\sigma)^{\alpha-1}\int_{\mathbb{R}^N}|u(\sigma,x)|^p\,dx\,d\sigma\,ds\nonumber\\
&=& \int_0^t\left(\int_{\sigma}^t(t-s)^{(N-5)/2}(s-\sigma)^{\alpha-1}\,ds\right)\int_{\mathbb{R}^N}|u(\sigma,x)|^p\,dx\,d\sigma \nonumber\\
&=&C \int_0^t(t-\sigma)^{(N-5)/2+\alpha}\left(\int_{|x|< r+\sigma}|u(\sigma,x)|^p\,dx\right)\,d\sigma\nonumber\\
&\geq&C (r+t)^{-N(p-1)}\int_0^t(t-s)^{(N-5)/2+\alpha}\left(\int_{\mathbb{R}^N}|u(s,x)|\,dx\right)^p\,ds.
\end{eqnarray}
Using H\"older's inequality, we have
\begin{equation}
\int_0^t(t-s)^{(N-5)/2}\left(\int_{\mathbb{R}^N}|u(s,x)|\,dx\right)\,ds\leq C (r+t)^{(N-3)(p-1)/(2p)-\alpha/p}\left(\int_0^t(t-s)^{(N-5)/2+\alpha}\left(\int_{\mathbb{R}^N}|u(s,x)|\,dx\right)^p\,ds\right)^{1/p},
\end{equation}
where we have used the fact that $\gamma>1-(p-1)(N-3)/2.$ Then
\begin{equation}\label{estimationB}
\int_0^t(t-s)^{(N-5)/2+\alpha}\left(\int_{\mathbb{R}^N}|u(s,x)|\,dx\right)^p\,ds\geq \frac{C|F(t)|^p}{(r+t)^{(N-3)(p-1)/2-\alpha}}.
\end{equation}
Therefore, combining $(\ref{estimationA})$ and $(\ref{estimationB}),$ we obtain
\begin{equation}\label{estimationC}
\ddot{F}(t)\geq \frac{C|F(t)|^p}{(r+t)^{(N-3)(p-1)/2+N(p-1)-\alpha}}.
\end{equation}
On the other hand, by repeating the same calculation in \cite[Section 5 p. 391]{Sideris}, we have
\begin{equation}\label{estimationD}
\int_{t-r}^t(t-s)^{(N-5)/2}\int_{|x|>t}u(s,x)\,dx\,ds\geq \int_{t-r}^t(t-s)^{(N-5)/2}\int_{|x|>t}v(s,x)\,dx\,ds,\qquad t\;\mbox{large},
\end{equation}
where $v$  is the solution of the homogeneous equation
$$
\left\{
\begin{array}{ll}
\displaystyle{v_{tt}-\Delta v=0}&\displaystyle{(x,t)\in\mathbb{R}^N\times(0,T),}\\
\displaystyle{v(x,0)=u_0(x),\;v_t(x,0)=u_1(x)}&\displaystyle{x\in\mathbb{R}^N,}\\
\end{array}
\right.
$$ Using H\"older's inequality and the compact support of $u$ on the left of $(\ref{estimationD}),$ one has
\begin{eqnarray}\label{estimationE}
\int_{t-r}^t(t-s)^{(N-5)/2}\int_{|x|>t}u(s,x)\,dx\,ds&\leq& \int_0^t(t-s)^{(N-5)/2}\int_{t<|x|<r+s}|u(s,x)|\,dx\,ds\nonumber\\
&\leq&C (r+t)^{(N-3)(p-1)/(2p)-\alpha/p}\left(\int_0^t(t-s)^{(N-5)/2+\alpha}\left(\int_{t<|x|<r+t}|u(s,x)|\,dx\right)^p\,ds\right)^{1/p}\nonumber\\
&\leq&C (r+t)^{(N-3)(p-1)/(2p)-\alpha/p}(r+t)^{(N-1)(p-1)/p}\left(\ddot{F}(t)\right)^{1/p}.
\end{eqnarray}
Next, to estimate the right sided of $(\ref{estimationD}),$ it follows from $(\ref{averages1})$ and \cite[Lemma 6]{Sideris} that
\begin{equation}\label{estimationF}
 \int_{t-r}^t(t-s)^{(N-5)/2}\int_{|x|>t}v(s,x)\,dx\,ds\geq C(r+t)^{(N-1)/2},\qquad t\;\mbox{large}.
\end{equation}
Hence, $(\ref{estimationD}),$ $(\ref{estimationE})$ and $(\ref{estimationF})$ imply
$$\ddot{F}(t)\geq (r+t)^{N-1-(N-1)p/2-(N-3)(p-1)/2+\alpha},$$
which leads after two integrations, for large $t,$ that
\begin{equation}\label{estimationG}
F(t)\geq (r+t)^{N+1-(N-1)p/2-(N-3)(p-1)/2+\alpha},
\end{equation}
where we have used the fact that $N+1-(N-1)p/2-(N-3)(p-1)/2+\alpha>1.$
Finally, making use of \cite[Lemma 4]{Sideris}, it follows from $(\ref{estimationC})$ and $(\ref{estimationG})$ that $T_{\max}<\infty,$ provided $p<p_2.$ $\hfill\square$\\
\end{pf}

\begin{theorem}\label{Blow-up-theorem3}$(\gamma\leq{(N-2)}/{N}\;\mbox{and $N\geq3$})$\\
Let $N\geq3$ and $p>1$ satisfies $(\ref{conditionoverp0}).$ Assume that $0<\gamma\leq {(N-2)}/{N}$ and $(u_0,u_1)\in H^1(\mathbb{R}^N)\times L^2(\mathbb{R}^N)$ such that
$$
\int_{\mathbb{R}^N}u_0>0,\quad \int_{\mathbb{R}^N}u_1>0.
$$
 If $p\leq {1}/{\gamma},$ then the solution of $(\ref{1++})$
blows up in finite time.\\
\end{theorem}
\begin{pf} Let $u$ be a global weak solution of $(\ref{1++}).$ Our argument is the same one of Theorem \ref{blowuptheorem++}.
 So we have two cases:\\

\noindent $\bullet$ \underline{The case $p\;< (1 / \gamma )$}: We repeat
the same argument as in the case $p<p_1,$ introduced in Theorem \ref{blowuptheorem++}, by choosing the following
function $\tilde{\varphi}(x,t)=\left(\varphi_1(x)\right)^\ell\varphi_2(t)$
where $\varphi_1(x)=\Phi\left({|x|}/{R}\right),$
$\varphi_2(t)=\left(1- {t}/{T}\right)^\eta_+$, $\ell,\eta\gg1$ and $R\in(0,T)$ large
enough such that when $T\rightarrow\infty$ we don't have $R\rightarrow\infty$ at the same time, with the same function $\Phi.$ We then obtain
\begin{eqnarray}\label{e1004++}
&{}&\int_{\mathcal{C}_T} |u|^p\;\tilde{\varphi}
\;+\;C\;T^{-\alpha}\int_{\mathcal{C}}\left(\varphi_1(x)\right)^\ell u_1(x)\;+\;C\;T^{-\alpha-1}\int_{\mathcal{C}}\left(\varphi_1(x)\right)^\ell u_0(x)\nonumber\\
&{}&\leq\;
\int_{\mathcal{C}_T}|u|\;\tilde{\varphi}^{1/p}\tilde{\varphi}^{-1/p}\left(\varphi_1(x)\right)^\ell
\left|D^{2+\alpha}_{t|T}\varphi_2(t)\right|\nonumber\\
&{}&\quad+\;C\int_{\mathcal{C}_T}|u|\;\tilde{\varphi}^{1/p}\tilde{\varphi}^{-1/p}(\varphi_1(x))^{\ell-1}
\left|\Delta_x\varphi_1(x)\;D^\alpha_{t|T}\varphi_2(t)\right|\nonumber\\
&{}&\quad+\;C\int_{\mathcal{C}_T}|u|\;\tilde{\varphi}^{1/p}\tilde{\varphi}^{-1/p}(\varphi_1(x))^{\ell-2}
|\nabla\varphi_1(x)|^2\left|D^\alpha_{t|T}\varphi_2(t)\right|,
\end{eqnarray}
\noindent where
$$\mathcal{C}_T:=[0,T]\times\mathcal{C}:=[0,T]\times\left\{x\in{\mathbb{R}}^N\hspace{2
mm};\hspace{2mm}|x|\leq 2R\right\},\;\;\int_{\mathcal{C}_T}=\int_{\mathcal{C}_T}\,dx\,dt,\;\;\int_{\mathcal{C}}=\int_{\mathcal{C}}\,dx.$$
 \noindent Now, by Young's
inequality $(\ref{2.16++}),$ with the same $a$ and $b$ as above and using $(\ref{positivity1++}),$ we get
\begin{eqnarray*}
\int_{\mathcal{C}_T}
|u|^p\;\tilde{\varphi}&\leq&C\int_{\mathcal{C}_T}
\left(\varphi_1(x)\right)^\ell \left(\varphi_2(t)\right)^{-\frac{1}{p-1}}
\left|D^{2+\alpha}_{t|T}\varphi_2(t)\right|^{\tilde{p}}\\
&{}&+\;C\int_{\mathcal{C}_T}\left(\varphi_1(x)\right)^{\ell-\tilde{p}}
\left(\varphi_2(t)\right)^{-\frac{1}{p-1}}
\left|\Delta_x\varphi_1(x)D^\alpha_{t|T}\varphi_2(t)\right|^{\tilde{p}}\\
&{}& +\;C\int_{\mathcal{C}_T}\left(\varphi_1(x)\right)^{\ell-2\tilde{p}}
\left(\varphi_2(t)\right)^{-\frac{1}{p-1}}
|\nabla\varphi_1(x)|^2\left|D^\alpha_{t|T}\varphi_2(t)\right|^{\tilde{p}}.
\end{eqnarray*}
\noindent Then, the new variables $\xi=R^{-1}x,$ $\tau=T^{-1}t$ and formulas
$(\ref{2.14++})$ and $(\ref{1000++})$ allow us to obtain
\begin{equation}\label{conditionABC}
\int_{\mathcal{C}_T}|u(x,t)|^p\;\tilde{\varphi}(x,t) \leq
\;C\;T^{1-(2+\alpha)\tilde{p}}\;R^N\;+\;C\;T^{1-\alpha
\tilde{p}}\;R^{N-2\tilde{p}}.
\end{equation}
Taking the limit as $T\rightarrow\infty,$ we infer, as $p<\frac{1}{\gamma}$ $(\Longleftrightarrow\;\;1-\alpha\tilde{p}<0),$ that
$$\int_0^\infty\int_{\mathcal{C}}|u(x,t)|^p\left(\varphi_1(x)\right)^\ell\,dx\,dt=0.$$
\noindent Finally, by taking $R\rightarrow\infty,$ we get a contradiction.\\

\noindent $\bullet$ \underline{The case $p\;= (1 / \gamma )$}: Here, we take the same test function in the last case. So, from $(\ref{conditionABC}),$ we obtain
$$\int_{\mathcal{C}_T}|u(x,t)|^p\;\tilde{\varphi}(x,t) \leq
\;C\;T^{-2\tilde{p}}\;R^N\;+\;C\;R^{N-2\tilde{p}}.$$
Taking the limit as $T\rightarrow\infty,$ we infer
$$\int_0^\infty\int_{\mathcal{C}}|u(x,t)|^p\left(\varphi_1(x)\right)^\ell\,dx\,dt\leq C\;R^{N-2\tilde{p}}.$$
Now, as the conditions ${(N-2)}/{N}<\gamma<1$ and $p={1}/{\gamma}$ imply that $N-2\tilde{p}<0,$ therefore, after passing to the limit as $R\rightarrow\infty,$
we conclude that
$$\int_0^\infty\int_{\mathbb{R^N}}|u(x,t)|^p\,dx\,dt =0;$$
contradiction and our result is established. $\hfill\square$  \\
\end{pf}




\section*{Acknowledgements}
The first named author was supported by the Lebanese National Council of Scientific Research (CNRS). He would like to thank very much the second named author for his hospitality during his stay in Pisa. The third  author was supported by the Italian National Council of Scientific Research (project PRIN No.
2008BLM8BB )
entitled: "Analisi nello spazio delle fasi per E.D.P."
\bibliographystyle{elsarticle-num}



\end{document}